\documentclass[pdflatex,sn-mathphys-num]{sn-jnl}

\usepackage{graphicx}%
\usepackage{multirow}%
\usepackage{amsmath, amssymb, amsthm, mathtools,bm}
\usepackage{mathrsfs}%
\usepackage[title]{appendix}%
\usepackage{xcolor}%
\usepackage{textcomp}%
\usepackage{manyfoot}%
\usepackage{algorithm, algorithmicx, algpseudocode}%
\usepackage{geometry}
\usepackage{hyperref}
\usepackage{comment}
\usepackage{lmodern}
\usepackage{tabularx, booktabs}
\usepackage{enumitem}

\newcommand{\vect}[1]{\bm{#1}}        
\newcommand{\mat}[1]{\mathbf{#1}}     
\newcommand{\op}[1]{\mathcal{#1}}     

\theoremstyle{thmstyleone}%
\newtheorem{theorem}{Theorem}
\newtheorem{corollary}{Corollary}
\newtheorem{lemma}{Lemma}
\newtheorem{Assumption}{Assumption}
\newtheorem{proposition}{Proposition}%

\theoremstyle{thmstyletwo}%
\newtheorem{remark}{Remark}%

\theoremstyle{thmstylethree}%

\begin{document}

\title[PG-Flow: Deterministic implicit policy gradients for geometric product-form queueing networks]{PG-Flow: Deterministic implicit policy gradients for geometric product-form queueing networks \thanks{This manuscript is an extended preprint version of the under review version}}



\author[1]{\fnm{Youssef} \sur{Ait El Mahjoub}}\email{youssef[dot]ait-el-mahjoub[at]efrei[dot]fr}

\affil[1]{\orgdiv{Efrei Research Lab}, \orgname{Universit\'e Paris-Panth\'eon-Assas}, \orgaddress{\street{30-32 Av. de la R\'epublique}, \city{Villejuif}, \postcode{94800}, \country{France}}}


\abstract{
Product-form queueing networks (PFQNs) admit steady-state distributions that factorize into local terms, and in many classical PFQNs—including Jackson, BCMP, G-networks, and Energy Packet Networks—these marginals are geometric and parametrized by local flow variables satisfying balance equations. While this structure yields closed-form expressions for key performance metrics, its use for deterministic steady-state optimization remains limited. We introduce \emph{PG--Flow}, a deterministic policy-gradient framework that differentiates through the steady-state flow fixed-point equations, providing exact gradients via implicit differentiation and a local adjoint system while avoiding trajectory sampling and Poisson equations. We establish global convergence under structural assumptions (affine flow operators and convex local costs), and show that acyclic networks admit linear-time computation of both flows and gradients. Numerical experiments on routing control in Jackson networks and energy-arrival control in Energy Packet Networks demonstrate that PG--Flow provides a principled and computationally efficient approach to deterministic steady-state optimization in geometric product-form networks.
}


\maketitle

\begingroup
\renewcommand\thefootnote{\fnsymbol{footnote}}
\footnotetext[1]{This manuscript is an extended preprint version of the under review document.}
\endgroup

\section{Introduction}

Stochastic networks provide a fundamental framework for modeling the dynamics
of complex systems such as communication infrastructures, manufacturing lines,
and energy-aware sensor networks. Among these models, \emph{product-form
queueing networks} (PFQNs) occupy a central place due to their analytical
tractability: under classical quasi-reversibility and independence conditions,
their steady-state probability distribution factorizes into a product of local
terms~\cite{jackson1957networks,BCMP1975,Gelenbe1991}. This property enables
efficient computation of steady-state metrics such as throughput, occupancy, or
delay, and has led to extensive work on performance evaluation and sensitivity
analysis in both open and closed networks~\cite{Gelenbe1993,
boucherie2011fundamental}. More recently, the class of systems
with product-form steady-state distributions has expanded through
model-specific quasi-reversibility arguments, as in pass-and-swap
queues~\cite{comte2021passswap}. Early contributions primarily focused on
characterizing and evaluating steady-state behavior under fixed system
parameters~\cite{jackson1957networks,BCMP1975,Gelenbe1991,Gelenbe1993,
GELENBE2002179,Gelenbe2012,Gelenbe2016}.

As applications grew in scale and complexity, attention shifted from analysis
to \emph{control} and \emph{optimization} of stochastic networks. In such
settings, routing probabilities, service rates, or energy-allocation policies
may be adjusted to improve long-run performance. This raises the challenge of
optimizing steady-state metrics that depend implicitly on control parameters
through the global balance equations of the network.
In parallel, stochastic control and reinforcement learning have developed two major families of methods for optimizing decision-making in controlled Markov processes: value-based approaches and policy-based approaches~\cite{puterman1994mdp,bertsekas1999nonlinear}.
Value-based techniques, such as classical policy iteration, rely on repeated solutions of Poisson or Bellman equations to evaluate a fixed policy and then perform an improvement step. While theoretically sound, these methods become computationally intractable in large state spaces, motivating efforts to exploit structural decompositions of the underlying MDP.
A first line of work addresses scalability through hierarchical MDPs, in which the decision problem is decomposed into subtasks; notable examples include MAXQ value-function decomposition~\cite{dietterich2000} and hierarchical abstraction frameworks~\cite{barto2003}. A complementary direction exploits factored or structured MDPs, where the transition and reward functions admit a graphical or algebraic decomposition, enabling low-dimensional value-function representations~\cite{koller1999}. In the same spirit of leveraging structure for computational efficiency, several works have focused on structured decision processes arising in energy-aware systems, such as slot-based or superstate-decomposable MDPs~\cite{ait2025slot,ait2024finding,ait2025efficient}, which allow scalable optimization in domains where standard DP becomes prohibitive.
Despite these advances, these techniques rely on discrete state–action spaces and depend heavily on dynamic-programming operators. When the underlying dynamics or action spaces are continuous, these methods typically require approximation schemes, such as diffusion-based limits, linear approximation architectures, or neural-network parameterizations.
To overcome such limitations, policy-based methods—most notably policy-gradient and actor–critic algorithms—estimate gradients directly from stochastic trajectories~\cite{sutton1999policy,konda2000actor}. More recently, implicit-differentiation techniques~\cite{domke2012implicit,gould2019deep} have enabled the computation of exact gradients through equilibrium constraints, bridging reinforcement learning with modern optimization theory.

Closer to product-form networks, several works have studied control based on
interaction with the steady-state regime. Sanders et al.~\cite{sanders2012online}
proposed an online stochastic-approximation scheme for adjusting the parameters
of a product-form network using empirical steady-state observations. Their
approach exploits the structural decomposition of PFQNs and relies on noisy
gradient estimates obtained from time-window sampling. More recently, Comte and collaborators introduced the \textit{Score-Aware Gradient Estimator (SAGE)}~\cite{comte2025sage}, a variance-reduced stochastic policy-gradient method based on Monte-Carlo sampling, applicable when the stationary distribution belongs to an exponential family (as in certain product-form systems). They establish local convergence in probability under the assumption that a suitable Lyapunov function exists. A complementary line of work focuses on controlling stochastic networks through
Lyapunov-drift techniques, as developed by Neely~\cite{neely2010stochastic}.
This framework yields robust online algorithms—such as MaxWeight or
drift–plus–penalty scheduling—that optimize long-run time averages without
computing steady-state distributions. While highly effective for dynamic
control, these approaches do not exploit product-form structures and do not
address gradient-based steady-state optimization.

Despite these advances, existing methods lack \emph{deterministic} gradient
formulas that explicitly leverage the internal structure of product-form
networks. In many classical PFQNs—including Jackson, BCMP, G-networks, and
Energy Packet Networks—the steady-state distribution factorizes into geometric
marginals parametrized by local flow variables. These flows satisfy a system of
local balance equations linking arrival rates, service rates, activation
probabilities, and routing. Although such equations are central to product-form
derivations, they have been rarely used as the computational basis for policy
optimization. Yet, in these models, performance metrics can be expressed solely
in terms of local flow variables, enabling implicit differentiation through the
associated equilibrium equations.

To address these limitations, we introduce \emph{PG-Flow}, a framework for
\emph{deterministic policy-gradient optimization} tailored to product-form
networks with geometric steady-state marginals. PG-Flow replaces global Poisson
equations and sampling-based estimators with a system of local steady-state
flow equations that implicitly characterize the long-run behavior of the
network. Through implicit differentiation and the solution of a local adjoint
system, PG-Flow computes exact gradients with quasi-linear complexity in the
number of queues in certain structured cases. In the general setting, PG--Flow operates as a deterministic implicit-gradient
method built upon the steady-state flow equations of a product-form network.  
Its computational and convergence properties are governed by the structural
features of these flow equations: affine operators yield favorable convexity,
while acyclic topologies induce triangular dependencies enabling efficient
solvers.  
These structural regimes illustrate how analytical queueing-network models with
product-form steady-state distributions can support scalable steady-state
optimization.  
We conclude the introduction by summarizing our main contributions.
\begin{itemize}
    \item We introduce a flow-based implicit-differentiation framework for geometric product-form networks: In PFQNs whose steady-state distribution factorizes into geometric marginals, 
    several classical performance metrics such as queue length, delay, utilization, 
    or leakage, admit natural closed-form expressions in terms of local steady-state 
    flows.  
    Building on this structure, we develop an implicit-differentiation approach and 
    present \emph{PG--Flow}, a deterministic policy-gradient algorithm that computes 
    exact steady-state gradients through the local equilibrium equations.

    \item We establish global convergence under structural assumptions and identify regimes of linear-time complexity: We prove that, under standard conditions (openness of the network, 
    monotonicity of the flow operator, convexity of local costs, and in particular 
    \emph{affinity} of the steady-state flow equations), the steady-state performance 
    objective becomes globally convex, and PG--Flow converges to the unique global 
    minimizer.  
    Furthermore, in \emph{acyclic} networks the flow operator becomes triangular, 
    making both the fixed-point computation and the adjoint solution reducible to 
    back-substitution, and yielding an overall \emph{linear-time} per iteration  algorithm.

    \item We demonstrate PG--Flow on two representative control problems.
    We apply our framework to (i) a routing-control problem in an open Jackson 
    network, and (ii) an energy-arrival control problem in an Energy Packet Network 
    (EPN).   These examples illustrate how PG--Flow exploits product-form structure to 
    compute deterministic policy gradients efficiently in both classical and 
    energy-aware queueing systems.
\end{itemize}

Table~\ref{tab:related} summarizes the conceptual distinctions between existing
approaches and the present work.

\begin{table}[hbtp]
\centering
\caption{Comparison of approaches for steady-state policy optimization in stochastic networks.}
\label{tab:related}
\begin{tabularx}{\linewidth}{lXXXX}
\toprule
\textbf{Method} &
\textbf{Evaluation} &
\textbf{Gradient} &
\textbf{Use of product-form} &
\textbf{Guarantees} \\
\midrule

Ref.~\cite{puterman1994mdp} &
Poisson / Bellman equation for policy evaluation &
Not gradient-based; exact DP updates &
None &
Global optimality; monotone policy improvement \\ \hline

Ref.~\cite{neely2010stochastic} &
Lyapunov drift-plus-penalty based on queue backlogs &
No steady-state gradient; drift / subgradient-type control &
None &
Queue stability and asymptotic optimality of time averages \\ \hline

Ref.~\cite{sanders2012online} &
Empirical steady-state sampling over time windows &
Stochastic gradient (Robbins--Monro) &
Exploits product-form steady-state distribution (incl.\ normalization) &
Asymptotic convergence under mixing assumptions \\ \hline

Ref.~\cite{tsai2023} &
Analytical and QNA-based performance formulas &
No analytic gradient; simulated annealing over routes &
No explicit product-form (fork--join FIFO / infinite-server network) &
Heuristic optimization; no formal gradient convergence result \\ \hline

Ref.~\cite{comte2025sage} &
Monte-Carlo sampling of steady-state Markov chain &
Score-based policy gradient (likelihood ratio, SAGE) &
Requires product-form / Exponential-family steady-state &
Stochastic (local) convergence w.h.p. under local Lyapunov-stability assumptions \\ \hline

PG-Flow (this work) &
Local steady-state flow equations with geometric marginals &
Deterministic gradient (implicit differentiation + adjoint) &
Explicitly uses product-form / Geometric-family steady-state &
Local convergence in general; global under affine flows and (H1)--(H4); faster global convergence in acyclic networks \\
\bottomrule
\end{tabularx}
\end{table}

The remainder of the paper is organized as follows.
Section~\ref{sec:problem} introduces the flow-based representation of
geometric product-form networks and the associated steady-state fixed-point
equations.  
Section~\ref{sec:pg-flow} presents the PG--Flow algorithm, including the
implicit-gradient formulation and the adjoint system.  
Section~\ref{sec:global-convergence} establishes local and global convergence
properties under structural assumptions, and Section~\ref{sec:complexity}
analyzes the computational complexity and identifies the acyclicity structure yielding
linear-time implementations.  
Section~\ref{sec:numerical} illustrates PG--Flow on routing control in an open
Jackson network and on energy-arrival control in an Energy Packet Network.
Section~\ref{sec:conclusion} concludes the paper and discusses perspectives for
future research.

\section{Steady-state flow problem}
\label{sec:problem}

We consider an open stochastic network composed of $d$ interacting queues and
$p$ control parameters.  
Each queue $i\in\{1,\dots,d\}$ is associated with a nonnegative
\emph{local steady-state flow variable} $\phi_i\in\mathbb{R}_+$, representing the
steady-state intensity with which tasks (jobs, packets, or energy units,
depending on the model) enter and leave that queue. Here, the term \textit{steady-state flow} refers to the effective arrival/output rate of each queue in steady state.  
Under a control vector $\vect{\theta}\in\mathbb{R}^p$, this equilibrium flow is
governed by a local operator $\mathcal{G}_i(\vect{\phi};\vect{\theta})$, and
collecting all components yields the global vector
\begin{equation}
    \vect{\phi}=(\phi_1,\ldots,\phi_d)^\top \in \mathbb{R}_+^d .
\end{equation}
\paragraph{Product-form networks with geometric marginals.}
In this work we focus on networks whose steady-state joint probability distribution admits a
product-form representation with geometric marginals.  
Under a fixed policy $\vect{\theta}$, the steady-state probability of observing a
state $x=(x_1,\ldots,x_d)$ factorizes as
\begin{equation}
\label{eq:product-form}
\Pi_{\vect{\theta}}(x)
   = \prod_{i=1}^d
      \bigl(1-\rho_i^\star(\vect{\theta})\bigr)\,
      \rho_i^\star(\vect{\theta})^{\,x_i},
\end{equation} 
where $\rho_i^\star(\vect{\theta})\in[0,1)$ denotes the load of queue~$i$ in 
steady state.
This geometric form arises in several well-established product-form models, 
including Jackson networks, BCMP queues, G-networks, and Energy Packet 
Networks (EPNs).  
In each case, the parameter $\rho_i^\star(\vect{\theta})$ is typically a deterministic function of the local steady-state flow 
 $\phi_i^\star(\vect{\theta})$ and model-specific parameters as 
service rate or energy parameters of node~$i$. Because the geometric marginal depends essentially on these flows, many standard
steady-state performance measures (mean queue length, mean delay, loss rate, or energy efficiency, etc.) can be expressed directly in terms of
the steady-state flows. 

\paragraph{Steady-state flows via fixed-point relations.}
Although the network admits a product-form distribution, explicit computation of
the normalizing constant $C$ is not required for steady-state analysis.  
Instead, the network dynamics enforce a system of fixed-point equations
characterizing the steady-state flows.  
For each node $i$, the equilibrium relation takes the form
\begin{equation}
\label{eq:fixed-point-system}
\phi_i^\star
   = \mathcal{G}_i(\vect{\phi}^\star;\vect{\theta}),
\qquad i=1,\ldots,d.
\end{equation}
The operator $\mathcal{G}_i$ incorporates the model-specific mechanisms that
govern the flow through queue $i$. A compact notation is obtained by defining the global operator
\begin{equation}
\mathcal{G}(\vect{\phi};\vect{\theta})
 := \bigl(
      \mathcal{G}_1(\vect{\phi};\vect{\theta}),
      \ldots,
      \mathcal{G}_d(\vect{\phi};\vect{\theta})
    \bigr)^\top,
\end{equation}
so that the equilibrium flows satisfy the vector fixed-point equation
\begin{equation}
\label{eq:fixed-point-G}
\vect{\phi}^\star(\vect{\theta})
   = \mathcal{G}\!\left(
        \vect{\phi}^\star(\vect{\theta});
        \,\vect{\theta}
     \right).
\end{equation}

Solving~\eqref{eq:fixed-point-G} yields the steady-state flow vector
$\vect{\phi}^\star(\vect{\theta})$, from which the performance measures relevant
to the optimization problem can be evaluated.  
The fixed-point representation therefore plays a fundamental role in our
framework, serving as the sole interface between the controlled parameters and
the steady-state behavior of the network. For instance, in a Jackson network 
with $M/M/1$ queues, and a control on routing probabilities $P(\vect{\theta})$, the steady-state
flow satisfies
\begin{equation} 
\phi_i(\vect{\theta})
    = \lambda_i^{\mathrm{ext}}
      + \sum_{j=1}^d P_{j,i}(\vect{\theta}) \,\phi_j(\vect{\theta}),
\end{equation}
and the mean queue length is
$m_i(\phi_i^\star)=\phi_i^\star/(\mu_i-\phi_i^\star)$. While in G-networks \cite{Gelenbe1991,Gelenbe1993,GELENBE2002179,fourneauYaem2017} flows incorporate  the effects of triggered
customer, and further in Energy Packet Networks \cite{Gelenbe2012,Gelenbe2016,YHJM20}, $\phi_i^\star$ incorporates both data flow and energy flow, which
jointly determine the steady-state occupancy of each node.

\paragraph{Steady-state objective.}
Because steady-state performance decomposes into local contributions, the global
objective (i.e. reward) can be written as
\begin{equation}
\label{eq:objective}
J(\vect{\theta})
   = \sum_{i=1}^{d} w_i \, r_i\!\bigl(\phi_i^\star(\vect{\theta}), \vect{\theta}\bigr)
   = \mathcal{F}\!\bigl(\vect{\phi}^\star(\vect{\theta}),\,\vect{\theta}\bigr),
\end{equation}
where $\mathcal{F}:\mathbb{R}^d\times\mathbb{R}^p\to\mathbb{R}$ is continuously
differentiable and $w_i\geq 0$ are importance weights.  
All dependence on steady-state performance therefore flows through the
steady-state vector $\vect{\phi}^\star(\vect{\theta})$, which makes the fixed-point
relation~\eqref{eq:fixed-point-G} the fundamental structural element for the proposed deterministic gradient-based optimization method.

Our aim is to determine a control vector $\vect{\theta}^\star$ that minimizes
the steady-state objective~\eqref{eq:objective}.  
In a gradient-based approach, the direction
$\nabla_{\vect{\theta}} J(\vect{\theta})$ provides the information required to
update $\vect{\theta}$ and improve system performance
\cite{boyd2004convex, bertsekas1999nonlinear, puterman1994mdp}.  
However, the mapping
$\vect{\theta}\mapsto\vect{\phi}^\star(\vect{\theta})$ is defined \emph{implicitly}
through the nonlinear fixed-point system~\eqref{eq:fixed-point-G}.  
Consequently, computing $\nabla_{\vect{\theta}}J(\vect{\theta})$ requires
handling the implicit dependence of $\vect{\phi}^\star$ on the control
parameters, which cannot be written in closed form.

In the next section, we exploit this implicit structure to derive an
\emph{implicit policy-gradient formula}, which leads naturally to the proposed
\textbf{PG--Flow} algorithm.  For clarity, we also specify the dimension
conventions that will be used throughout the derivations.

\begin{table}[!h]
\centering
\caption{Notation and dimensions used throughout the paper.}
\renewcommand{\arraystretch}{1.4}
\begin{tabular}{c c l}
\hline
\textbf{Symbol} & \textbf{Dimensions / Domain} & \textbf{Description} \\
\hline
$\vect{\phi}$ & $\mathbb{R}_+^{d\times 1}$ & Vector of steady-state flow variables \\

$\vect{\phi}^\star(\vect{\theta})$ & $\mathbb{R}_+^{d\times 1}$ & Fixed point satisfying $\vect{\phi}^\star = \op{G}(\vect{\phi}^\star;\vect{\theta})$ \\

$\vect{\theta}$ & $\mathbb{R}^{p\times 1}$ & Policy parameter vector (controls) \\

$\op{G}(\vect{\phi},\vect{\theta})$ & $\mathbb{R}_+^{d\times 1}$ & Global flow operator (mapping flows to flows) \\

$\partial_{\vect{\phi}} \op{G}$ & $\mathbb{R}^{d\times d}$ & Jacobian of $\op{G}$ w.r.t.\ flows \\

$\partial_{\vect{\theta}} \op{G}$ & $\mathbb{R}^{d\times p}$ & Jacobian of $\op{G}$ w.r.t.\ parameters \\[3pt]

$\op{F}(\vect{\phi},\vect{\theta})$ & $\mathbb{R}$ & Global performance objective \\

$\partial_{\vect{\phi}} \op{F}$ & $\mathbb{R}^{1\times d}$ & Row gradient of $\op{F}$ w.r.t.\ flows \\

$\partial_{\vect{\theta}} \op{F}$ & $\mathbb{R}^{1\times p}$ & Row gradient of $\op{F}$ w.r.t.\ parameters \\

$\dfrac{d \vect{\phi}^\star}{d\vect{\theta}}$ & $\mathbb{R}^{d\times p}$ & Sensitivity of steady-state flows \\[3pt]

$y$ & $\mathbb{R}^{d\times 1}$ & Adjoint variable solving $(I - \partial_{\vect{\phi}} \op{G})^\top y = (\partial_{\vect{\phi}} \op{F})^\top$ \\

$\nabla_{\vect{\theta}} J$ & $\mathbb{R}^{1\times p}$ & Row gradient of objective \\
\hline
\end{tabular}
\end{table}

\section{PG-Flow algorithm}
\label{sec:pg-flow}

\begin{Assumption}
\label{assump1}
For the remainder of this work, we make the following assumptions, needed for the derivation of the implicit gradient:  
\begin{enumerate}[label=(\roman*)]
    \item For any admissible parameter vector $\vect{\theta}$, the flow solution $\vect{\phi}^\star(\vect{\theta})$
    exists. This assumption will be shown to hold for open networks (see Section \ref{sec:global-convergence}).
    \item The global operator $\op{G}(\vect{\phi}; \vect{\theta})$ is $C^1$ in both arguments.
    Its Jacobians $\partial_{\vect{\phi}}\op{G}\in\mathbb{R}^{d\times d}$
    and $\partial_{\vect{\theta}}\op{G}\in\mathbb{R}^{d\times p}$
    exist and are continuous.
    \item  $\op{F}$ is continuously differentiable.
\end{enumerate}
\end{Assumption}

The previous section introduced the steady-state flow equations
$\vect{\phi}^\star(\vect{\theta})$
and the corresponding performance objective
$J(\vect{\theta})$.
Optimizing this objective requires the gradient
$\nabla_{\vect{\theta}} J(\vect{\theta})$.
Because $\vect{\phi}^\star(\vect{\theta})$ depends implicitly on~$\vect{\theta}$
through the fixed-point system, the derivative must be propagated
through the equilibrium mapping~$\op{G}$.
We now derive an exact expression for this gradient and show how it leads
to an efficient algorithmic formulation.

Differentiating~\eqref{eq:objective} using the chain rule yields
\begin{equation}
\label{eq:chain}
\nabla_{\vect{\theta}} J
\;=\;
\partial_{\vect{\theta}}\op{F}
\;+\;
\partial_{\vect{\phi}}\op{F}\;
\frac{d\vect{\phi}^\star}{d\vect{\theta}},
\end{equation}
where $\nabla_{\vect{\theta}} J \in \mathbb{R}^{1\times p}$, 
$\partial_{\vect{\theta}}\op{F} \in \mathbb{R}^{1\times p}$,
$\partial_{\vect{\phi}}\op{F} \in \mathbb{R}^{1\times d}$,
and
$\frac{d\vect{\phi}^\star}{d\vect{\theta}} \in \mathbb{R}^{d\times p}$.

\begin{proof}
By total differentiation of
$J(\vect{\theta}) = \op{F}(\vect{\phi}^\star(\vect{\theta}), \vect{\theta})$,
we have
\begin{equation*}
dJ
\;=\;
\partial_{\vect{\phi}}\op{F}\; d\vect{\phi}^\star
\;+\;
\partial_{\vect{\theta}}\op{F}\; d\vect{\theta},
\end{equation*}
Since
\begin{equation*}
d\vect{\phi}^\star
\;=\;
\Big(\frac{d\vect{\phi}^\star}{d\vect{\theta}}\Big)\, d\vect{\theta},
\end{equation*}
we obtain
\begin{equation*}
dJ
\;=\;
\Big[
  \partial_{\vect{\theta}}\op{F}
  +
  \partial_{\vect{\phi}}\op{F}\,
  \frac{d\vect{\phi}^\star}{d\vect{\theta}}
\Big]\,
d\vect{\theta}.
\end{equation*}
Identifying the coefficient of $d\vect{\theta}$ yields
$\nabla_{\vect{\theta}} J = \frac{dJ}{d\vect{\theta}}$,
which proves~\eqref{eq:chain}.
\end{proof}

\paragraph{Closed-form differentiation vs. numerical recomputation.}
By \emph{closed-form} we mean an explicit expression for $\vect{\phi}^\star(\vect{\theta})$ that can be differentiated.
In the considered networks, this is typically intractable. A naive alternative is to \emph{recompute the fixed point}
for perturbed parameters and approximate derivatives by \textbf{finite differences}
\begin{equation}
\label{eq:fd}
\frac{\partial J}{\partial \theta_j}
\;\approx\;
\frac{J(\vect{\theta}+\varepsilon\,\vect{e}_j) - J(\vect{\theta}-\varepsilon\,\vect{e}_j)}{2\varepsilon},
\qquad j=1,\ldots,p.
\end{equation}
which requires recomputing the steady-state flows
$\vect{\phi}^\star(\vect{\theta}\pm \varepsilon \vect{e}_j)$ for \emph{every}
parameter direction.
Since each evaluation of $J(\vect{\theta}\pm\varepsilon \vect{e}_j)$ requires a full
solution of the fixed-point system~\eqref{eq:fixed-point-G}, finite differences incur
a cost proportional to $2p$ steady-state-flow solves per iteration, which rapidly
dominates the computational complexity. Implicit differentiation avoids this computational increase: the fixed-point system is solved
\emph{once} per iteration to obtain $\vect{\phi}^\star(\vect{\theta})$, after which the
local Jacobians $\partial_{\vect{\phi}}\op{G}$ and $\partial_{\vect{\theta}}\op{G}$ are evaluated around that solution.
Even when these Jacobians must be estimated numerically (e.g., by finite differencing
the map $\op{G}$ itself), one only evaluates $\op{G}$—not the full fixed-point
equation—thus avoiding repeated global solves.
This distinction is crucial: implicit differentiation separates the expensive operation
(solve $\vect{\phi}^\star$ once) from cheaper local sensitivity computations, whereas
finite differences repeatedly recompute the global equilibrium.
Section~\ref{sec:complexity} gives a full comparison of the resulting computational
costs.

\paragraph{Implicit differentiation.}
Instead of recomputing the fixed point for each perturbation of $\vect{\theta}$,
we use the \emph{implicit function theorem} (\cite[Prop.~A.25]{bertsekas1999nonlinear}). That is, 
differentiating \eqref{eq:fixed-point-G} locally yields
\begin{equation}
\label{eq:dphi-normal}
\frac{d\vect{\phi}^\star}{d\vect{\theta}}
\;=\;
\big(\mat{I}-\partial_{\vect{\phi}}\op{G}\big)^{-1}\,\partial_{\vect{\theta}}\op{G}.
\end{equation}
Substituting \eqref{eq:dphi-normal} into \eqref{eq:chain} yields the \emph{implicit gradient}:
\begin{equation}
\label{eq:grad-normal}
\nabla_{\vect{\theta}} J
\;=\;
\partial_{\vect{\theta}}\op{F}
\;+\;
\partial_{\vect{\phi}}\op{F}\,
\big(\mat{I}-\partial_{\vect{\phi}}\op{G}\big)^{-1}\,
\partial_{\vect{\theta}}\op{G}.
\end{equation}

\paragraph{Adjoint.}
Equation \eqref{eq:grad-normal} involves the inverse of 
$\mat{I}-\partial_{\vect{\phi}}\op{G}$. 
Rather than forming this inverse explicitly which is numerically less stable, in practice we introduce an adjoint variable 
$\vect{y}\in\mathbb{R}^{d\times 1}$ defined as the solution of
\begin{equation}
\label{eq:adjoint}
\big(\mat{I}-\partial_{\vect{\phi}}\op{G}\big)^{\!\top} \vect{y}
\;=\;
\big(\partial_{\vect{\phi}}\op{F}\big)^{\!\top}.
\end{equation}
This construction corresponds to applying the inverse 
$(\mat{I}-\partial_{\vect{\phi}}\op{G})^{-1}$ indirectly via the solution of the
linear system \eqref{eq:adjoint}. Hence, by right-multiplying \eqref{eq:grad-normal} by the transpose of the associated linear system. The transpose appears because the gradients are  \emph{row}  vectors while the adjoint system acts on 
\emph{column} vector $\vect{y}$. Substituting the definition of $\vect{y}$ into 
\eqref{eq:grad-normal} gives the \textit{compact implicit gradient}:
\begin{equation}
\label{eq:grad-adjoint}
\nabla_{\vect{\theta}} J
\;=\;
\partial_{\vect{\theta}}\op{F}
\;+\;
\vect{y}^{\!\top}\,\partial_{\vect{\theta}}\op{G}.
\end{equation}
Thus, the adjoint $\vect{y}$ evaluates
$\partial_{\vect{\phi}}\op{F}(\mat{I}-\partial_{\vect{\phi}}\op{G})^{-1}$
without explicit inversion, by solving a single local linear system of size $d$.

\paragraph{The PG-Flow algorithm}
By summarizing previous steps, we propose the following deterministic implicit policy-gradient algorithm: 
\begin{algorithm}[!ht]
\caption{PG-Flow: Policy Optimization via steady-state Flow Equations \label{algo:pg-flow}}
\begin{algorithmic}[1]
\Require Initial parameters $\vect{\theta}$
\While{not converged}
  \State Solve steady-state flows: $\vect{\phi}^\star=\op{G}(\vect{\phi}^\star;\vect{\theta})$.
  \State Evaluate the average cost: $J(\vect{\theta})=\sum_i w_i\,r_i(\phi_i^\star; \vect{\theta})$.
  \State Compute local derivatives: $\partial_{\vect{\phi}}\op{F}$, $\partial_{\vect{\theta}}\op{F}$, $\mat{G}_\phi{=}\partial_{\vect{\phi}}\op{G}$, $\mat{G}_\theta{=}\partial_{\vect{\theta}}\op{G}$.
  \State Solve the adjoint: $(\mat{I}-\mat{G}_\phi)^\top \, \vect{y} = (\partial_{\vect{\phi}}\op{F})^\top$.
  \State Global gradient: $\nabla_{\vect{\theta}} J = \partial_{\vect{\theta}}\op{F} + \vect{y}^\top\mat{G}_\theta$.
  \State Update: $\vect{\theta} \leftarrow \vect{\theta} - \eta\, \nabla_{\vect{\theta}} J$.
\EndWhile
\end{algorithmic}
\end{algorithm}

Sections \ref{sec:problem}–\ref{sec:pg-flow} only require Assumption \ref{assump1}. In Section \ref{sec:global-convergence} we introduce additional structural conditions (H1–H4) under which the PG-Flow iteration converges globally to a steady-state optimum.

\section{Global Convergence of PG-Flow}
\label{sec:global-convergence}

The PG-Flow algorithm is a gradient-based optimization scheme.
As such, it may in general converge only to a \emph{local} minimum of the steady-state objective $J(\vect{\theta})$.
To guarantee convergence to a \emph{global} optimum, it is therefore necessary to identify
sufficient conditions under which $J(\vect{\theta})$ is \emph{convex} on its control domain.

Let $\Phi \coloneqq \mathbb{R}_+^d$ denote the nonnegative orthant, representing all admissible steady-state flow
vectors, and let $\mathcal{U}\subset\mathbb{R}^p$ be a convex set of feasible control parameters (for instance a
simplex when a global resource budget must be split across nodes, as in the EPN example). We focus on a broad class
of stochastic networks---including Jackson, BCMP, G-networks, and energy-based product-form models---whose
steady-state flows $\vect{\phi}^\star(\vect{\theta})\in\Phi$ satisfy an \emph{affine} fixed-point relation:

\begin{equation}
\label{eq:flow-fixed}
\vect{\phi}^\star(\vect{\theta})
\;=\;
\op{G}(\vect{\phi}^\star(\vect{\theta}), \vect{\theta}),
\qquad
\op{G}(\vect{\phi}, \vect{\theta})
= A(\vect{\theta})\,\vect{\phi} + b(\vect{\theta}),
\end{equation}
Here:
\begin{itemize}
    \item $\vect{\theta}\in\op{U}$ denotes the tunable control parameters (routing probabilities, service rates, energy-injection rates, etc.).
    \item $A(\vect{\theta})\in\mathbb{R}_+^{d\times d}$ encodes the internal propagation of steady-state flow.  
    Its entry $A_{ij}(\vect{\theta})$ represents the contribution of the steady-state flow at node $j$ to that of node $i$.  
    In classical Jackson networks, $A(\vect{\theta})$ reduces exactly to the transpose routing matrix $P^\top$.  
    In BCMP networks, its block structure represents class-transition probabilities.  
    In G-networks, $A(\vect{\theta})$ may incorporate positive and negative triggered movements.  
    In Energy Packet Networks (EPNs) (for instance \cite{YHJM20}), $A(\vect{\theta})$ captures the steady-state DP flow propagation, which again reduces to $P^\top$ for the Data Packet queue.
    \item $b(\vect{\theta})\in\mathbb{R}_+^d$ collects all exogenous or self-generated contributions to the flow that are 
    independent of the internal interactions represented by $A(\vect{\theta})$.  
    In Jackson networks, $b(\vect{\theta})$ equals the external arrival rates.  
    In BCMP or G-networks, it may include class-dependent arrivals or triggered events as signals.  
    In EPNs, it includes the exogenous injections of data or energy packets.
\end{itemize}
Thus, the pair $(A(\vect{\theta}),b(\vect{\theta}))$ fully describes the internal flow propagation and the external inputs of the
steady-state regime. This general affine representation encompasses a wide family of product-form networks and forms
the basis for analyzing convexity and deriving global convergence guarantees for PG--Flow.

The global performance objective is expressed as a sum of local reward functions \eqref{eq:objective} where each $r_i$ is related to queue $i \in \{1,\ldots,d\}$.
Our goal is to determine the optimal control vector
\begin{equation}
\vect{\theta}^\star \in \arg\min_{\vect{\theta} \in \mathcal{U}} J(\vect{\theta}).
\end{equation}

In the following, we introduce structural assumptions (H1)–(H4) under which
$J(\vect{\theta})$ is convex on $\mathcal{U}$,
ensuring that PG-Flow converges to the global minimizer $\vect{\theta}^\star$.

\begin{itemize}
\item[(H1)] \textbf{Convex and compact}: 
The control set $\mathcal U$ is compact and convex.

\item[(H2)] \textbf{Openness}: 
For all $\vect{\theta} \in \mathcal{U}$, the network is open: from
every queue there exists a directed path to a departure node, so that no closed
communicating class of queues can retain flow indefinitely.

\item[(H3)] \textbf{Monotonicity}: 
The mapping $\vect{\theta} \mapsto \vect{\phi}^\star(\vect{\theta})$
is continuous on $\mathcal{U}$ and coordinate-wise non-decreasing.

\item[(H4)] \textbf{Local rewards}: 
Each local reward \( r_i(\phi_i,\vect{\theta}) \) is convex in $(\phi_i,\vect{\theta})$
and non-decreasing in $\phi_i$.
For example, \( r_i(\phi_i) = \phi_i/(\mu_i - \phi_i) \) on $[0,\mu_i[$.
\end{itemize}

\begin{lemma}[Existence of steady-state flows]
\label{lem:open-neumann}
If the network is open, then the steady-state flow vector $\vect{\phi}^\star$
exists and is unique.

\begin{proof}
Openness (H2) guarantees that no unit of flow can circulate indefinitely in the
network. This property is ensured under either one of the following 
structural conditions on the nonnegative interaction matrix $A(\vect{\theta})$:

\begin{itemize}
\item[(i)] \emph{Strictly sub-stochastic case:}
there exists $\kappa<1$ such that
$\sum_{j=1}^d A_{i,j}(\vect{\theta}) \le \kappa$ for all rows $i$.
Then $\|A(\vect{\theta})\|_\infty \le \kappa < 1$, so by Perron--Frobenius
$\rho(A(\vect{\theta})) < 1$.

\item[(ii)] \emph{Feed-forward (acyclic) case:}
$A(\vect{\theta})$ can be permuted into a strictly triangular form with zero
diagonal (no self-loops and no cycles). Triangularity implies that all
eigenvalues of $A(\vect{\theta})$ are zero, hence $\rho(A(\vect{\theta})) = 0$.
\end{itemize}

In both cases we obtain $\rho(A(\vect{\theta})) < 1$, which ensures that
$\mat{I} - A(\vect{\theta})$ is invertible and that the Neumann series converges:
\begin{equation*}
(\mat{I} - A(\vect{\theta}))^{-1}
= \sum_{n=0}^{\infty} A(\vect{\theta})^n.
\end{equation*}
Therefore the fixed-point equation
$\vect{\phi}^\star(\vect{\theta})
= A(\vect{\theta})\, \vect{\phi}^\star(\vect{\theta})
   + b(\vect{\theta})$
admits the unique solution
\begin{equation*}
\vect{\phi}^\star(\vect{\theta})
= (\mat{I} - A(\vect{\theta}))^{-1} b(\vect{\theta})
= \sum_{n=0}^{\infty} A(\vect{\theta})^n\, b(\vect{\theta}).\qedhere
\end{equation*}
\end{proof}
\end{lemma}

\begin{corollary} [Existence of adjoint solution]
\label{cor:adjoint} 
The adjoint system \eqref{eq:adjoint} admits a unique solution.
\begin{proof}
From the affine form in \eqref{eq:flow-fixed}, 
$\op{G}(\vect{\phi};\vect{\theta}) = A(\vect{\theta})\,\vect{\phi} + b(\vect{\theta})$,
so that
\[
\mat{G}_\phi(\vect{\theta}) = \frac{\partial \op{G}(\phi;\vect{\theta})}{\partial\vect{\phi}}
= A(\vect{\theta}),
\qquad \forall\, \vect{\phi},\vect{\theta}.
\]
and by Lemma~\ref{lem:open-neumann}, the matrix $I-A(\vect{\theta})$ is invertible for all $\vect{\theta}$, and therefore so is its transpose
$(I-A(\vect{\theta}))^\top$.
As a result, the adjoint system
\begin{equation*}
\bigl(I - \mat{G}_{\vect{\phi}}(\vect{\theta})\bigr)^\top\, \vect{y}
= (I-A(\vect{\theta}))^\top \vect{y} = \bigl(\partial_{\vect{\phi}} \op{F}(\vect{\phi}^\star(\vect{\theta}),\vect{\theta})\bigr)^\top,
\end{equation*}
admits a unique solution $y(\vect{\theta})\in\mathbb R^d$ for every $\vect{\theta}\in\mathcal U$.
\end{proof}
\end{corollary}

\begin{lemma}[Regularity and boundedness]
\label{lem:regul}
Let $\op{U}\subset\mathbb{R}^p$ be compact. 
Assume that $A(\vect{\theta})$ is continuously differentiable in~$\vect{\theta}$ on a neighborhood of $\op{U}$.
By Lemma~\ref{lem:open-neumann}, $I-A(\vect{\theta})$ is invertible for all $\vect{\theta}\in \op{U}$.
Then:
\begin{itemize}
    \item[(a)] the mapping $\vect{\theta}\mapsto (I-A(\vect{\theta}))^{-1}$ is continuously differentiable on $\op{U}$;
    \item[(b)] the fundamental matrices are uniformly bounded:
    \[
    \sup_{\vect{\theta}\in \op{U}}\|(I-A(\vect{\theta}))^{-1}\| < \infty.
    \]
\end{itemize}

\begin{proof}
(a)  
Since $A(\vect{\theta})$ is $C^1$, the map $\vect{\theta}\mapsto I-A(\vect{\theta})$ is $C^1$.  
Lemma~\ref{lem:open-neumann} ensures that $I-A(\vect{\theta})$ is invertible for every $\vect{\theta}\in \op{U}$.  
For an invertible matrix $M$, the identity
\begin{equation*}
\frac{d}{d\theta}(M^{-1}(\vect{\theta})) = - M^{-1}(\vect{\theta})\, M'(\vect{\theta})\, M^{-1}(\vect{\theta})
\end{equation*}
holds whenever $M'(\vect{\theta})$ exists.  
Applying this with $M(\vect{\theta})=I-A(\vect{\theta})$, we see that 
$(I-A(\vect{\theta}))^{-1}$ has a continuous derivative on $\op{U}$, hence is $C^1$.

(b)
The map $\theta\mapsto (I-A(\vect{\theta}))^{-1}$ is continuous on the compact set $\op{U}$,  
so its norm attains a finite maximum.
\end{proof}
\end{lemma}

\begin{lemma}[Smoothness of the steady-state flow]
\label{lem:smooth}
Let $\vect{b}(\vect{\theta})$ be continuously differentiable on a neighborhood of $\op{U}$.
For each $\vect{\theta}\in \op{U}$, the steady-state flow is given by
\begin{equation}
\vect{\phi}^{\star}(\vect{\theta})
= \big(I - A(\vect{\theta})\big)^{-1} \vect{b}(\vect{\theta}) .
\end{equation}
Under the assumptions of Lemma~\ref{lem:regul}, the mapping 
$\vect{\theta}\mapsto \vect{\phi}^{\star}(\vect{\theta})$ is $C^1$ on $\op{U}$, 
and its Jacobian is uniformly bounded on~$\op{U}$:
\begin{equation}
\sup_{\vect{\theta}\in \op{U}}
\big\|\, d\vect{\phi}^{\star}(\vect{\theta}) \,\big\| < \infty .
\end{equation}
\begin{proof}
From Lemma~\ref{lem:regul}, the inverse matrix $(I-A(\vect{\theta}))^{-1}$ is $C^1$ in $\vect{\theta}$
and its norm is uniformly bounded on $\op{U}$.
Since $\vect{b}(\vect{\theta})$ is also $C^1$, the product
$\vect{\phi}^{\star}(\vect{\theta})
= (I-A(\vect{\theta}))^{-1}\, \vect{b}(\vect{\theta})$
is $C^1$ as a composition and product of $C^1$ mappings.

Using the standard product rule for matrix–vector functions,
\begin{equation*}
d\vect{\phi}^{\star}(\vect{\theta})
=
d\!\left((I-A(\vect{\theta}))^{-1}\right)\,\vect{b}(\vect{\theta})
\;+\;
(I-A(\vect{\theta}))^{-1}\, d\vect{b}(\vect{\theta}) .
\end{equation*}

Each term on the right-hand side is continuous on the compact set $\op{U}$:
\begin{itemize}
\item $(I-A(\vect{\theta}))^{-1}$ is continuous and uniformly bounded (Lemma~\ref{lem:regul}),
\item $d\!\left((I-A(\vect{\theta}))^{-1}\right)$ is continuous on $\op{U}$ (Lemma~\ref{lem:regul}(a)),
\item $\vect{b}(\vect{\theta})$ and $d\vect{b}(\vect{\theta})$ are continuous by assumption.
\end{itemize}

A continuous function on a compact set is bounded.
Hence both terms above are uniformly bounded on~$\op{U}$, and so is
$d\vect{\phi}^{\star}(\vect{\theta})$.
\end{proof}
\end{lemma}

\begin{proposition}[Lipschitz continuity of the steady-state gradient]
\label{prop:lips}
Assume each local reward $r_i(\phi_i,\vect{\theta})$ is twice continuously differentiable 
on its domain. Define 
$J(\vect{\theta}) = \sum_{i=1}^d w_i\, r_i(\phi_i^\star(\vect{\theta}),\vect{\theta})$, $w_i \geq 0$.
Under Lemma~\ref{lem:regul} and Lemma~\ref{lem:smooth}, the mapping $J$ is twice continuously differentiable on $\op{U}$.
Moreover, its Hessian is uniformly bounded on $\op{U}$, i.e.
\begin{equation}
\label{eq:hess-bound}
L = \sup_{\vect{\theta}\in \op{U}}\|\nabla^2 J(\vect{\theta})\| < \infty .
\end{equation}
Consequently, the gradient $\nabla J$ is Lipschitz continuous on $\op{U}$:
\begin{equation}
\|\nabla J(\vect{\theta}) - \nabla J(\vect{\theta}')\| 
\le L \|\vect{\theta}-\vect{\theta}'\| ,\qquad \forall \vect{\theta},\vect{\theta}'\in \op{U} ,
\end{equation}

for some finite constant $L$.

\begin{proof}
By Lemma~\ref{lem:smooth}, $\vect{\phi}^\star(\vect{\theta})$ is $C^1$ with uniformly bounded
Jacobian on~$\op{U}$. 
Since each $r_i(\phi_i,\vect{\theta})$ is $C^2$, the composite mapping
$\vect{\theta} \;\longmapsto\; r_i\big(\phi_i^\star(\vect{\theta}),\,\vect{\theta}\big)$
is also $C^{2}$. Multiplying by a nonnegative constant $w_i$ preserves
$C^2$-smoothness. Therefore the finite sum
$J(\vect{\theta})=\sum_{i=1}^d w_i\, r_i\bigl(\phi_i^\star(\vect{\theta}),\vect{\theta}\bigr)$
is $C^{2}$ on $\op{U}$. Moreover, Each entry of $\nabla^2 J(\vect{\theta})$ is obtained from compositions 
and products of continuous functions (the partial derivatives of $r_i$, which are 
continuous by the $C^2$ assumption, and the derivative of $\vect{\phi}^\star$, 
which is continuous by Lemma~\ref{lem:smooth}). 
Hence each entry of $\nabla^2 J(\vect{\theta})$ is continuous on the compact set $\op{U}$
and its norm achieves a finite maximum \eqref{eq:hess-bound}.

A standard result in smooth optimization then yields the Lipschitz continuity of the
gradient: if the Hessian of a $C^2$ function is upper bounded by~$L$ on a convex set,
then its gradient is $L$-Lipschitz on that set  
(\cite[Eq.~(9.13)]{boyd2004convex}).  
Therefore,
\begin{equation*}
J(\vect{\theta}') 
\;\le\; 
J(\vect{\theta}) 
+ \nabla J(\vect{\theta})^T(\vect{\theta}'-\vect{\theta})
+ \frac{L}{2}\|\vect{\theta}'-\vect{\theta}\|^2,
\end{equation*}
and 
\begin{equation*}
\|\nabla J(\vect{\theta}) - \nabla J(\vect{\theta}')\|
\;\le\;
L\,\|\vect{\theta}-\vect{\theta}'\|,
\qquad
\forall \vect{\theta},\vect{\theta}'\in U .
\end{equation*}
Thus $\nabla J$ is $L$-Lipschitz on $\op{U}$.
\end{proof}
\end{proposition}

\begin{theorem}[Convexity of $J$ and global convergence of PG-Flow]
\label{th:global-PGFlow}
Under assumptions {\rm(H1)}–{\rm(H4)}, the steady-state flow mapping 
$\vect{\theta}\mapsto \vect{\phi}^\star(\vect{\theta})$ defined by 
$\vect{\phi}^\star = A(\vect{\theta})\,\vect{\phi}^\star + b(\vect{\theta})$
is well defined on $\mathcal U$. Moreover, with $w_i\geq0$ the steady-state objective
$J(\vect{\theta})
= \sum_{i=1}^d w_i r_i\!\big(\phi_i^\star(\vect{\theta}),\,\vect{\theta}\big)$
is convex on $\mathcal U$ and admits at least one global minimizer
$\vect{\theta}^\star\in\arg\min_{\mathcal U}J$.
In addition, $J$ is continuously differentiable on $\mathcal U$ and its gradient is 
Lipschitz continuous on $\mathcal U$. Therefore, for any initialization $\vect{\theta}_0\in\mathcal U$, the PG-Flow iteration
\begin{equation}
\vect{\theta}_{k+1} 
= \vect{\theta}_k - \eta_k \,\nabla J(\vect{\theta}_k)
\end{equation}
with constant step-size $\eta_k\equiv\eta\in(0,2/L)$, or with an Armijo
backtracking rule, produces a non-increasing sequence $J(\vect{\theta}_k)$ and 
converges to a global minimizer of $J$ on $\mathcal U$.

\begin{proof}
\textit{Step 1 (Existence and uniqueness of steady-state flows).}
By Lemma~\ref{lem:open-neumann} and Assumption (H2), $I-A(\vect{\theta})$ is invertible for all
$\vect{\theta}\in\mathcal U$.
Hence the fixed-point equation admits the unique solution
$\vect{\phi}^\star(\vect{\theta}) = (I-A(\vect{\theta}))^{-1}b(\vect{\theta})$.

\smallskip
\textit{Step 2 (Existence of the adjoint solution).}
By Corollary~\ref{cor:adjoint}, the adjoint system \eqref{eq:adjoint}
admits a unique solution $y(\vect{\theta})$ for all $\vect{\theta}\in\mathcal U$,
and the implicit gradient \eqref{eq:grad-adjoint}
$\nabla J(\vect{\theta})$
is therefore well defined.

\smallskip
\textit{Step 3 (Convexity of $J$).}
By Assumption (H3), $\vect{\theta}\mapsto \vect{\phi}^\star(\vect{\theta})$ is 
coordinatewise nondecreasing.
By Assumption (H4), each local reward $r_i(x,\vect{\theta})$ is convex in $(x,\vect{\theta})$ 
and nondecreasing in $x$.
Hence by composition with the monotone map $\phi_i^\star$, 
each function 
$\vect{\theta}\mapsto r_i(\phi_i^\star(\vect{\theta}),\vect{\theta})$ 
is convex.
Summing with $w_i \geq 0$ preserves convexity and yields the convexity of $J$ on $\mathcal U$.

\smallskip
\textit{Step 4 (Differentiability and Lipschitz continuity of the gradient).}
By Lemma~\ref{lem:smooth} the mapping $\vect{\theta}\mapsto \vect{\phi}^\star(\vect{\theta})$ is $C^1$ 
with uniformly bounded Jacobian on $\mathcal U$.
Since each $r_i$ is $C^2$, Proposition~\ref{prop:lips} ensures that $J$ is $C^2$ on $\op{U}$, with a uniformly bounded Hessian \eqref{eq:hess-bound}. By the quadratic upper bound \cite[Eq.\,(9.13)]{boyd2004convex}, $\nabla J$ is $L$-Lipschitz on $\mathcal U$.

\smallskip
\textit{Step 5 (Monotone descent and convergence).}
Since $\nabla J$ is $L$-Lipschitz, the "descent lemma" in 
\cite[Prop.~A.24]{bertsekas1999nonlinear}
implies that for all $x,y$ with $x,x+y\in\mathcal U$,
\begin{equation*}
\label{eq:descent-lemma}
J(x+y) \;\le\; 
J(x) + \nabla J(x)^{\!\top}y + \tfrac{L}{2}\|y\|^2.
\end{equation*}

Applying \eqref{eq:descent-lemma} with $x=\vect{\theta}_k$ and
$y=-\eta\,\nabla J(\vect{\theta}_k)$ yields 
\begin{align}
J(\vect{\theta}_{k+1})
&= J\big(\vect{\theta}_k - \eta\nabla J(\vect{\theta}_k)\big) \nonumber\\
&\le J(\vect{\theta}_k)
   +\big(-\eta\nabla J(\vect{\theta}_k)\big)^{\top}\nabla J(\vect{\theta}_k)
   +\frac{L}{2}\big\|-\eta\nabla J(\vect{\theta}_k)\big\|^{2}\nonumber\\
&= J(\vect{\theta}_k)
   -\eta\,\|\nabla J(\vect{\theta}_k)\|^{2}
   +\frac{L}{2}\,\eta^{2}\,\|\nabla J(\vect{\theta}_k)\|^{2}\nonumber\\
&= J(\vect{\theta}_k)
   -\Big(\eta-\tfrac{L}{2}\eta^{2}\Big)\,\|\nabla J(\vect{\theta}_k)\|^{2}.
\label{eq:PGFlow-descent}
\end{align}
Hence $\{J(\vect{\theta}_k)\}$ is non-increasing for any constant
step-size $\eta\in(0,2/L)$. Since $J$ is convex on $\mathcal U$ and has an $L$-Lipschitz continuous
gradient, the standard convergence theory of gradient descent
(see, e.g. \cite[Sec.~9.1]{boyd2004convex} or \cite{bertsekas1999nonlinear}) implies that
$\vect{\theta}_k$ converges to a global minimizer
$\vect{\theta}^\star\in\arg\min_{\vect{\theta}\in\mathcal U} J(\vect{\theta})$.
\end{proof}
\end{theorem}
\begin{remark}[Projected PG-Flow under physical constraints]
In many models the feasible control set $\mathcal{U}\subset\mathbb{R}^p$
is a compact convex subset reflecting physical limitations.  
Typical examples include:  

(i) constrained routing probabilities (section~\ref{subsec:jack});

(ii) admission or allocation rates constrained by capacity limits (section~\ref{subsec:EPN}); 

(iii) energy allocation vectors subject to a global budget constraint (section~\ref{subsec:EPN}),  
e.g. 
\begin{equation}
\mathcal{U}=\Bigl\{\vect{\theta}\in\mathbb{R}^p_+:
\sum_{i=1}^p \theta_i \le B_{\max}\Bigr\}.
\end{equation}

The unconstrained PG-Flow update
$\vect{\theta}_{k+1}
= \vect{\theta}_k - \eta \,\nabla J(\vect{\theta}_k)$
may leave $\mathcal{U}$.  
To maintain feasibility, the update is replaced by the 
\emph{Euclidean projected} iteration
\begin{equation}
\vect{\theta}_{k+1}
=\Pi_{\mathcal{U}}\!\bigl(\vect{\theta}_k-\eta\nabla J(\vect{\theta}_k)\bigr),
\qquad
\Pi_{\mathcal{U}}(x)=\arg\min_{u\in\mathcal{U}} \|u-x\|_2.
\end{equation}

Convexity of $\mathcal{U}$ ensures that $\Pi_{\mathcal{U}}$ is single-valued
and non-expansive 
($\|\Pi_{\mathcal{U}}(x)-\Pi_{\mathcal{U}}(y)\|\le\|x-y\|$).  
Thus the descent arguments in the proof of 
Theorem~\ref{th:global-PGFlow} remain valid for the projected scheme,
and global convergence continues to hold.
\end{remark}

\section{Computational Complexity Analysis \label{sec:complexity}}
The computation of the steady-state objective 
\( J(\vect{\theta}) \)
and its gradient depends critically on the method used to evaluate 
\(\nabla_{\vect{\theta}} J\).  
Because \(J\) depends on the steady-state flow vector 
\(\vect{\phi}^\star(\vect{\theta})\), any gradient method must address two 
fundamental operations:
(i) solving the steady-state fixed-point equation 
$\vect{\phi}^\star = \op{G}(\vect{\phi}^\star;\vect{\theta})$ and  
(ii) propagating derivatives through the implicit dependence of 
\(\vect{\phi}^\star\) on~\(\vect{\theta}\).  
Different approaches treat these two operations very differently, which leads to 
significant variations in computational cost. In this section, we compare five gradient-estimation methods relevant to 
steady-state optimization in product-form queueing networks:
finite differences applied to the full objective (FD–J),
Monte–Carlo policy gradients,
PG–Flow with numerical Jacobians,
PG–Flow with analytical Jacobians,
and finally PG–Flow in feed-forward networks.
Each method reflects a distinct philosophy regarding how the steady state 
interacts with the control parameters, and therefore provides a different trade-off 
between accuracy, and computational efficiency.

Evaluating the steady-state objective 
\(J(\vect{\theta}) = F(\vect{\phi}^\star(\vect{\theta}),\vect{\theta})\)
requires solving the fixed-point equation $\vect{\phi}^\star = \op{G}(\vect{\phi}^\star;\vect{\theta})$. In the dense case, the evaluation of a single application 
\(\op{G}(\vect{\phi},\vect{\theta})\) costs \(O(d^2)\), reflecting the internal 
dependencies among queues. While in a local-based case, each queue interacts only with a limited neighbors queues. 
Two classes of solvers may be used for obtaining \(\vect{\phi}^\star\):
(i) \textit{Direct linear solvers}, which reformulate the fixed point 
    as a linear system in the affine case and cost \(O(d^3)\). (ii) \textit{Iterative fixed-point solvers}, such as Picard iteration or 
    Anderson acceleration, whose per-iteration cost is \(O(d^2)\), and whose 
    total cost can be written as 
    $C_{\text{FP}} = O(K_{\text{fp}} d^2).$

In PG–Flow we favour iterative solvers, as they preserve sparsity, avoid matrix 
factorization, and lead to better scalability in large networks.

\subsection{Finite Differences Applied to the Objective (FD–J)}

A classical baseline approach consists of approximating 
\(\nabla_{\vect{\theta}}J\) through centered finite differences \eqref{eq:fd}:
each component \(\partial J/\partial \theta_j\) requires the computation of 
\(J(\vect{\theta}+\varepsilon e_j)\) and \(J(\vect{\theta}-\varepsilon e_j)\).
Thus every gradient evaluation requires \(2p\) steady-state solves:
\begin{equation}
\label{eq:O-FD-J}
\text{Cost}_{\text{FD–J}}
    = O(2p\, C_{\text{FP}}) 
    = O(2p\, K_{\text{fp}} d^2).
\end{equation}

Although straightforward, this approach is computationally expensive because the 
entire steady-state computation must be repeated for every coordinate direction.  
This becomes prohibitive for large-scale systems with many control parameters.

\subsection{Monte–Carlo Policy Gradient}

For comparison, Monte–Carlo policy-gradient methods avoid differentiating the 
steady state entirely. They rely instead on the identity
\begin{equation}
\nabla_{\vect{\theta}} J(\vect{\theta})
    = \mathbb{E}\!\left[
         \nabla_{\vect{\theta}} 
            \log \pi_{\vect{\theta}}(a|s)\, Q^{\pi}(s,a)
      \right],
\end{equation}
whose cost is dominated by the sampling of trajectories.  
If \(N_{\text{traj}}\) trajectories of horizon \(H\) are required, then
\begin{equation}
\label{eq:O-MC}
\text{Cost}_{\text{MC}} = O(N_{\text{traj}} H).
\end{equation}

While this estimator avoids manipulating the fixed-point equations, its variance 
grows with the mixing time of the queueing process.
This phenomenon is particularly problematic in heavy-traffic or highly 
congested networks, where correlations become strong and Monte-Carlo policy 
gradients suffer from slow and unstable convergence.
This motivates the development of deterministic, variance-free approaches such as 
PG–Flow.

\subsection{PG–Flow with Numerical Jacobians}

The PG–Flow algorithm \ref{algo:pg-flow} computes exact implicit gradients without perturbing the 
full objective, and without relying on samples.  
When analytical Jacobians are not available, the required matrices 
\(\partial_{\phi}\op{G}\) and \(\partial_{\vect{\theta}}\op{G}\) can still be recovered 
numerically by finite differences applied locally to \(\op{G}\), not to the full 
mapping \(\vect{\phi}^\star(\vect{\theta})\).  
This distinction is important: the fixed-point equation is solved only \emph{once} 
per iteration, after which all sensitivities are obtained from local evaluations 
of~\(\op{G}\).
Approximating the Jacobian \(\partial_{\phi}\op{G}\in\mathbb{R}^{d\times d}\)  requires computing \(d^2\)
partial derivatives, each obtained from one evaluation of \(\op{G}\). Since the
dense cost of evaluating \(\op{G}\) is \(O(d^2)\), we obtain
\[
\text{Cost}(\partial_{\phi}\op{G}) = O(d^4).
\]
Likewise, computing the \(dp\) entries of \(\partial_{\vect{\theta}}\op{G} \in\mathbb{R}^{d\times p}\)
costs
\[
\text{Cost}(\partial_{\vect{\theta}}\op{G}) = O(dp\, d^2).
\]
Once these Jacobians are available, the implicit gradient requires solving the 
adjoint linear system \eqref{eq:adjoint} whose cost is \(O(d^3)\) in the dense case.

Putting all contributions together, the dominant term is the numerical evaluation 
of the Jacobian \(\partial_{\phi}\op{G}\), leading to
\begin{equation}
\label{eq:O-numeric1}
\text{Cost}_{\text{PG--Flow,num}}
    = O(d^4 + dp\, d^2 + d^3)
    = O(d^4).\qquad (dense)
\end{equation}
In locally connected networks where evaluating \(\op{G}\) costs only \(O(d)\), the corresponding bound reduces to
\begin{equation}
\label{eq:O-numeric2}
\text{Cost}_{\text{PG--Flow,num}}
    = O(d^3 + dp\, d + d^3)
    = O(d^3).\qquad (locally\ connected)
\end{equation}

This cost remains significantly lower than FD–J for moderate \(p\), and benefits 
from the fact that only \(\op{G}\) (and not the full fixed point) is perturbed.

\subsection{PG–Flow with Analytical Jacobians}

In many product-form queueing networks (Jackson, BCMP, G-networks, 
Energy-Packet Networks), both \(\op{G}\) and the objective \(\op{F}\) admit 
closed-form derivatives.  
In such settings, each entry of the Jacobians 
\(\partial_{\phi}\op{G}\) and \(\partial_{\vect{\theta}}\op{G}\) can be computed in 
\(O(1)\), with sparsity determined by the local connectivity of the network.

Let \(n_z\) denote the total number of nonzero entries across these Jacobians.
In locally connected networks, \(n_z = O(d)\).  
The overall complexity is therefore
\begin{equation}
\text{Cost}_{\text{PG--Flow,analytic}}
    = O(n_z) + O(d^3),
\end{equation}
where the \(O(d^3)\) term corresponds to solving the adjoint system.
Hence, dense linear algebra dominates the iteration cost, and the Jacobian
construction becomes negligible in comparison.

\subsection{PG–Flow in Feed-Forward Networks}

A particularly favorable situation arises when the queueing network is 
\emph{feed-forward}, i.e.\ acyclic.  
In this case, the dependency structure of the flow operator is triangular, and 
each component \(\op{G}_i\) depends only on flows from its predecessors.
As a consequence, both the steady-state flows and their sensitivities propagate 
in a single pass. Hence, this structure has the following implications:
\begin{itemize}
    \item the forward computation of the steady-state flows 
    \(\vect{\phi}^\star\) requires a single pass \(O(d)\);
    \item the Jacobian \(\partial_{\vect{\phi}}\op{G}\) contains only \(O(d)\) nonzeros;
    \item the parameter Jacobian \(\partial_{\vect{\theta}}\op{G}\) is block-sparse and local;
    \item the adjoint linear system is triangular and solves in \(O(d)\) by back-substitution.
\end{itemize}

Combining these observations, the total cost of PG--Flow in a feed-forward 
network reduces to
\begin{equation}
\label{eq:O-analytic}
\text{Cost}_{\text{PG--Flow,analytic,acyclic}} = O(d),
\end{equation}
which represents the best achievable scaling: linear time per iteration.
This highlights the strong link between the queuing network structure and the computational 
efficiency of PG--Flow. We now illustrate the practical behavior of PG--Flow on representative numerical
examples in Section~\ref{sec:numerical}.

\begin{table}[h!]
\centering
\caption{Computational complexity of gradient-estimation methods for 
steady-state queueing networks.\label{tab:complexity}}
\renewcommand{\arraystretch}{1.3}
\begin{tabular}{lcc}
\toprule
\textbf{Method} &
\textbf{Gradient type} &
\textbf{Cost / iteration} \\
\midrule
Finite differences (FD--J) &
Perturb $J(\vect{\theta})$ &
$O(p\, d^2)$--$O(p\, d^3)$ \\[2pt]
PG-explicit (MC) &
Stochastic PG estimator &
$O(N_{\mathrm{traj}} H)$ \\[2pt]
PG--Flow (numeric, dense) &
Implicit + FD Jacobians &
$O(d^4)$ \\[2pt]
PG--Flow (numeric, locally connected) &
Implicit + FD Jacobians &
$O(d^3)$ \\[2pt]
PG--Flow (analytic) &
Implicit (closed form) &
$O(n_z) + O(d^3)$ \\[2pt]
PG--Flow (analytic + acyclic) &
Implicit (triangular/sparse) &
\textbf{$O(d)$} \\
\bottomrule
\end{tabular}
\end{table}

All PG--Flow variants that rely on exact implicit gradients (obtained either
analytically or from numerically exact Jacobians) inherit the global convergence
guarantees of Theorem~\ref{th:global-PGFlow} under assumptions (H1)--(H4).  
Acyclicity is not required for convergence: it only improves computational
complexity.  
Identifying broader structural classes of stochastic networks (e.g., hierarchical,
nearly triangular, decomposable or Nearly-Completely-Decomposable \cite{stewart1994markov}, with central cycles \cite{Rob90,ait2025slot,ait2025efficient}) for which both the steady-state
flow equations and the adjoint system admit sub-cubic or even linear solvers is an
interesting direction for future work.

We now illustrate the practical behavior of PG--Flow on representative numerical
examples in Section~\ref{sec:numerical}.

\section{Numerical analysis \label{sec:numerical}}


\subsection{Jackson network with controllable routing matrix}
\label{subsec:jack}

We consider a three-node feed-forward Jackson network with fixed external Poisson arrivals $\lambda_i^{\mathrm{ext}}$
and two internal routing parameters. Each node corresponds to an M/M/1 queue. Jobs arrive from outside only into queue~1 with rate
$\lambda_1^{\mathrm{ext}}=4$, while queues~2 and~3 receive no
external arrivals.
The service is exponential and rates are fixed as
\[
\mu = (\mu_1,\mu_2,\mu_3) = (6,5,7).
\]

Routing is controlled at nodes~1 and~2.
Upon service completion at queue~1, jobs are routed
\begin{itemize}
  \item to queue~2 with probability \(p(\theta_1)\),
  \item to queue~3 with probability \(1-p(\theta_1)\).
\end{itemize}
Upon service completion at queue~2, jobs are routed
\begin{itemize}
  \item to queue~3 with probability \(q(\theta_2)\),
  \item and leave the network with probability \(1-q(\theta_2)\).
\end{itemize}
Jobs leaving queue~3 always exit the system. The routing probability matrix is therefore expressed as:
\begin{equation*}
P(\vect{\theta})=
\begin{bmatrix}
0 & p(\theta_1) & 1-p(\theta_1)\\[3pt]
0 & 0 & q(\theta_2)\\[3pt]
0 & 0 & 0
\end{bmatrix}
\end{equation*}

The decision vector collects the two routing controls,
\[
\vect{\theta} = (\theta_1,\theta_2)^\top \in \mathcal U = \mathbb [0,1]^2.
\]

We use affine parameterizations
\[
p(\theta_1) = p_0 + \theta_1, \qquad
q(\theta_2) = q_0 + \theta_2,
\]
with nominal values \(p_0 = q_0 = 0\) and bounds on~\(\vect{\theta}\) chosen
so that $p(\theta_1),q(\theta_2)\in[0,1]$ for all
\(\vect{\theta} \in \mathcal U\).
In the numerical example below we take
\(\theta_1,\theta_2\in[0,1]\), so that \(p(\theta_1)=\theta_1\),
\(q(\theta_2)=\theta_2\), and \(p,q\in[0,1]\). Let \(\Lambda_i\) denote the steady-state incoming flow 
at queues~\(i=1,2,3\), and let
\(\phi = (\Lambda_1,\Lambda_2,\Lambda_3)^\top\).
The flow balance equations are
\[
\Lambda_1 = 4,\qquad
\Lambda_2 = p(\theta_1)\,\Lambda_1,\qquad
\Lambda_3 = \bigl(1-p(\theta_1)\bigr)\Lambda_1 + q(\theta_2)\Lambda_2.
\]
These can be written as a fixed-point relation
$\vect{\phi}^\star(\vect{\theta}) = \op{G}(\vect{\phi}^\star(\vect{\theta});\vect{\theta})
= A(\vect{\theta})\vect{\phi}^\star(\vect{\theta}) + b(\vect{\theta})$
 with $\vect{\phi}=[\phi_1, \phi_2, \phi_3]^T$, $A(\vect{\theta}) = P(\vect{\theta})^\top$ and 
\begin{equation}
\label{eq:num:fixed}
\op{G}(\vect{\phi};\vect{\theta}) =
\begin{bmatrix}
\lambda_1^{\mathrm{ext}}\\[3pt]
p(\theta_1)\,\phi_1\\[3pt]
\bigl(1-p(\theta_1)\bigr)\phi_1 + q(\theta_2)\phi_2
\end{bmatrix},
\qquad
A(\vect{\theta})=
\begin{bmatrix}
0 & 0 & 0\\[3pt]
p(\theta_1) & 0 & 0\\[3pt]
1-p(\theta_1) & q(\theta_2) & 0
\end{bmatrix},
\qquad
b(\vect{\theta}) = 
\begin{bmatrix}
\lambda_1^{\mathrm{ext}}\\0\\0
\end{bmatrix}.
\end{equation}

We aim to minimize the total mean queue length in the network. Let $\rho_i = \frac{\phi_i(\vect{\theta})}{\mu_i}$ be the load of queue i, then 
\begin{equation}
\label{eq:jackson-2param-cost}
J(\vect{\theta})
= \sum_{i=1}^3 \frac{\rho_i}{1 - \rho_i} = \sum_{i=1}^3 \frac{\Lambda_i^\star(\vect{\theta})}{\mu_i - \Lambda_i^\star(\vect{\theta})}.
\end{equation}
Detailed verifications of assumptions (H1)--(H4) and the explicit first PG--Flow
iteration are given in Appendix~\ref{annexe:A} (Jackson).

\begin{figure}[htbp]
\hspace{-1cm}
\includegraphics[width=15cm, height=6cm]{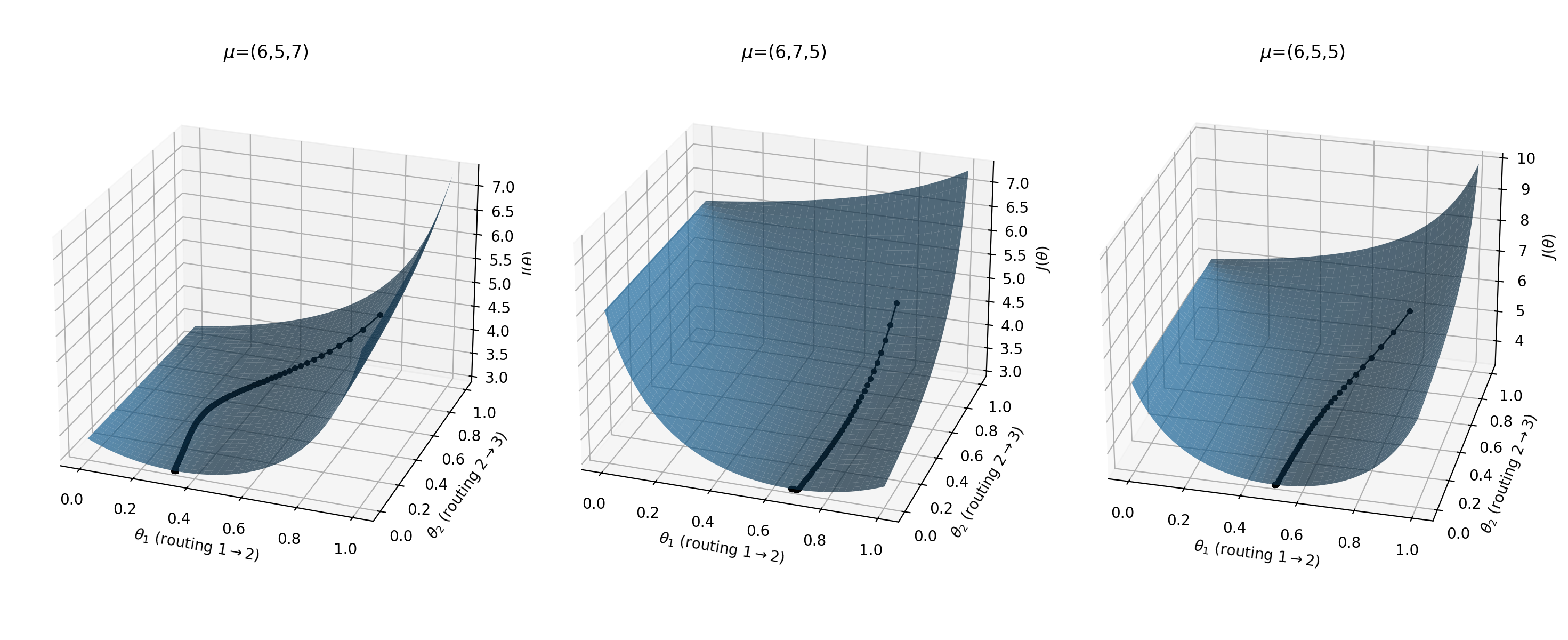}
\vspace{-1cm}
\caption{Cost surfaces $J(\theta_1,\theta_2)$ for the three Jackson configurations,
together with the PG-Flow trajectories (black curves).}
\label{fig:jackson-3scenarios}
\end{figure}

To illustrate how the PG-Flow algorithm behaves across different service-rate
configurations, Figure~\ref{fig:jackson-3scenarios} displays the cost surfaces 
$J(\vect{\theta})$ for the three scenarios considered (below), together with the 
algorithmic trajectories. Although the geometry of the landscape changes with $(\mu_2,\mu_3)$, the feasible
parameter domain remains convex in all cases, in agreement with our theoretical
assumptions.  

\paragraph{Scenario 1: server 3 is faster.}
In this experiment (the one detailed in Appendix \ref{annexe:A}), we considered a Jackson network with service
rates $\mu = (\mu_1,\mu_2,\mu_3) = (6,5,7),$
where server~3 is faster than server~2.  
The PG-Flow algorithm yields the optimal parameters
\[
\vect{\theta}^\star \approx (0.3313,\,0),\qquad
\vect{\phi}^\star \approx (4,\ 1.33,\ 2.67),
\qquad
J(\vect{\theta}^\star)\approx 2.979.
\]
Thus, the optimal policy routes only about one third of the input flow through queue~2,
while setting $q^\star=0$ to avoid sending customers from queue~2 to the slower 
queue~3. Since $\mu_3>\mu_2$, the optimal routing chooses predominantly the faster 
path $1\to 3$ and avoids the slower path $1\to2\to3$.

\paragraph{Scenario 2: server 2 becomes faster.}
To analyze how the optimal routing changes when service capacities are modified,
we consider the reversed situation
\[
\mu = (6,7,5),
\]
where server~2 is now faster than server~3.  
The PG-Flow algorithm converges to
\[
\vect{\theta}^\star \approx (0.6688,\,0),\qquad
\vect{\phi}^\star \approx (4,\ 2.68,\ 1.32),
\qquad
J(\vect{\theta}^\star)\approx 2.979.
\]
The optimal policy therefore \emph{switches side}: a large fraction of the 
incoming flow is now routed to queue~2, and only a small residual flow enters queue~3.  
As in Scenario~1, the algorithm finds $q^\star=0$, meaning that customers leave 
the system directly after queue~2, which is now the fastest internal server.  
Thus the optimized policy always favors routing through the fastest available server.

\paragraph{Scenario 3: balanced internal servers.}
Finally, we consider a configuration where the two internal servers have equal
service rates,
\[
\mu = (6,5,5).
\]
In this case PG-Flow converges to
\[
\vect{\theta}^\star \approx (0.5013,\,0),\qquad
\vect{\phi}^\star \approx (4,\ 2.01,\ 1.99),
\qquad
J(\vect{\theta}^\star)\approx 3.333.
\]
Since $\mu_2=\mu_3$, the optimal routing splits the flow almost equally between
queues~2 and~3.  However, because both servers are slower than in Scenarios~1--2,
the global performance is significantly worse: $J(\vect{\theta}^\star)$ is strictly larger.

\medskip
\noindent
\textbf{Overall comparison.}
Scenarios~1 and~2 achieve similar and significantly lower optimal costs
($J(\vect{\theta}^\star)\approx 2.98$), each by directing the majority of the flow toward
the fastest internal server.  
Scenario~3, where the internal servers are slower and balanced, yields a noticeably
higher cost ($J(\vect{\theta}^\star)\approx 3.33$).  
Hence the most favorable configuration is obtained when a single internal server
is sufficiently fast, allowing the optimizer to route most of the flow through it. Another observation is that 
in all three scenarios, the PG-Flow solution satisfies 
either \(p^\star=0\) or \(q^\star=0\).  
In these cases, the steady-state flows obey the conservation identity
\(\Lambda_2 + \Lambda_3 = \Lambda_1\).  
This matches the analytical expression
\[
\Lambda_1 = 4,\qquad 
\Lambda_2 = 4p,\qquad
\Lambda_3 = 4(1-p+pq)
\quad\Longrightarrow\quad
\Lambda_2+\Lambda_3 = 4 + 4pq,
\]
which reduces to \(\Lambda_2+\Lambda_3=\Lambda_1\) when
\(p=0\) or \(q=0\).  
Thus, the conservation pattern observed numerically is fully consistent
with the model’s flow equations.

\subsection{Energy Packet Network with controllable energy arrivals}
\label{subsec:EPN}

In this second application, we consider the Energy Packet Network (EPN) composed of $N$ interconnected  nodes presented in \cite{YHJM20}.
Each node $i \in \{1,\dots,N\}$ contains:
\begin{itemize}
    \item an infinite Data Packet (DP) queue with external Poisson arrivals of rate $\lambda_i^{\mathrm{ext}}$;
    \item an infinite Energy Packet (EP) queue receiving EPs at a \emph{controllable} Poisson rate 
    $\alpha_i = \alpha_i(\vect{\theta})$;
    \item a DP service is triggered by the EP queue: at each service event, one EP triggers the departure of a single DP from the DP queue (if a DP is available). 
          If no EPs are available, DPs wait in the DP queue until the next EP arrival. 
          After service, the EP is consumed and the DP is forwarded to another queue 
          according to the routing matrix $P$.
    \item an energy leakage mechanism in the EP queue where EPs disappear at exponential rate $\gamma_i$ (representing the physical energy leakage that may occur in IoT batteries).
\end{itemize}

We restrict ourselves to the Coxian case with one phase $K=1$ (i.e., exponential EP service times). The steady-state consumption rate of EPs at node $i$ is:
\begin{equation}
    \beta_i(\vect{\theta}) = \frac{\alpha_i(\vect{\theta})}{\gamma_i + \mu_i},
\end{equation}
where $\mu_i$ is the exponential EP service parameter. Thus, the EP flow depends affinely on the control vector $\vect{\theta}$.

\smallskip
\textbf{\textit{Routing of data packets.}}
Only DP flows are routed between nodes. Each DP served at node $i$ is forwarded to node $j$
with fixed probability $P_{i,j}$ or leaves the network with probability
$d_i = 1 - \sum_{j=1}^N P_{i,j}$. Routing probabilities are \emph{not} controlled in this experiment.
Instead, the control acts solely on the EP arrival rates $\alpha_i(\vect{\theta})$.
This highlights the versatility of PG--Flow, which can optimize arbitrary differentiable
steady-state parameters (routing, service rates, or energy input levels) without altering the
mathematical structure.

\smallskip
\textbf{\textit{Steady-state flow equations.}}
Let $\rho_i(\vect{\theta})$ denote the DP utilization at node $i$. Since the effective DP service rate
is $\mu_i \beta_i(\vect{\theta})$, the DP balance equations (special case $K=1$ of \cite{YHJM20})
are:
\begin{equation}\label{eq:EPN_rho}
\rho_i(\vect{\theta})
=
\frac{
\lambda_i^{\mathrm{ext}} 
+ 
\sum_{j=1}^N 
    \mu_j \rho_j(\vect{\theta}) \beta_j(\vect{\theta}) P_{j,i}
}{
\mu_i \beta_i(\vect{\theta})
},
\qquad i = 1,\dots,N.
\end{equation}

Multiplying \eqref{eq:EPN_rho} by $\mu_i\beta_i(\vect{\theta})$ yields:
\begin{equation}
\rho_i(\vect{\theta}) \mu_i \beta_i(\vect{\theta})
=
\lambda_i^{\mathrm{ext}}
+
\sum_{j=1}^N 
    \mu_j \rho_j(\vect{\theta})\beta_j(\vect{\theta})
    P_{j,i}.
\end{equation}

Define the DP flow:
\begin{equation}
    \phi_i(\vect{\theta})
    :=
    \rho_i(\vect{\theta}) \mu_i \beta_i(\vect{\theta}).
\end{equation}
Then the DP equations reduce to the classical Jackson flow equations:
\begin{equation}
\phi_i(\vect{\theta})
=
\lambda_i^{\mathrm{ext}} + \sum_{j=1}^N \phi_j(\vect{\theta}) P_{j,i}.
\end{equation}
Hence, DP flows satisfy the linear fixed-point system:
\begin{equation}
    \vect{\phi}^\star(\vect{\theta}) 
    = 
    A \vect{\phi}^\star(\vect{\theta}) + b,
\end{equation}
where $A = P^\top$ and $b = (\lambda_1^{\mathrm{ext}},\dots,\lambda_N^{\mathrm{ext}})^\top$.
Importantly, since routing is fixed, $A$ and $b$ do \emph{not} depend on $\vect{\theta}$. 
All dependence on $\vect{\theta}$ arises through the EP flows $\beta_i(\vect{\theta})$.

\smallskip
\textbf{\textit{Performance objective and QoS--Energy trade-off.}}
Two antagonistic performance metrics arise naturally in EPNs:

\begin{itemize}
    \item \textbf{DP queueing delay (Quality of Service)}:
    \[
        D(\vect{\theta})
        =
        \sum_{i=1}^N
        \frac{\rho_i(\vect{\theta})}{1 - \rho_i(\vect{\theta})}
        =
        \sum_{i=1}^N
        \frac{\phi_i(\vect{\theta})}
             {\mu_i\beta_i(\vect{\theta}) - \phi_i(\vect{\theta})},
    \]
    where the expression corresponds to an $M/M/1$ queue with service rate $\mu_i\beta_i(\vect{\theta})$.
    Increasing $\alpha_i(\vect{\theta})$ increases $\beta_i(\vect{\theta})$ and hence reduces delay.

    \item \textbf{EP leakage (energy waste)}:
    \[
        L(\vect{\theta})
        =
        \sum_{i=1}^N 
        \gamma_i \beta_i(\vect{\theta})
        =
        \sum_{i=1}^N 
        \gamma_i \frac{\alpha_i(\vect{\theta})}{\gamma_i+\mu_i}.
    \]
    Increasing $\alpha_i$ increases the energy waste.
\end{itemize}

These metrics are structurally conflicting: more incoming energy improves QoS but increases the
system’s energy consumption. We consider the scalarized objective:
\begin{equation}
    J(\vect{\theta})
    =
    w_1 D(\vect{\theta}) + w_2 L(\vect{\theta}),
    \qquad w_1,w_2 \ge 0,
\end{equation}
which expresses explicitly the QoS--Energy compromise.

Detailed verifications of assumptions (H1)--(H4) and the explicit first PG--Flow
iteration are given in Appendix~\ref{annexe:B} (EPN). Now, we present the following results. 
\begin{figure}[!ht]
    \centering
\includegraphics[width=10cm, height=5.5cm]{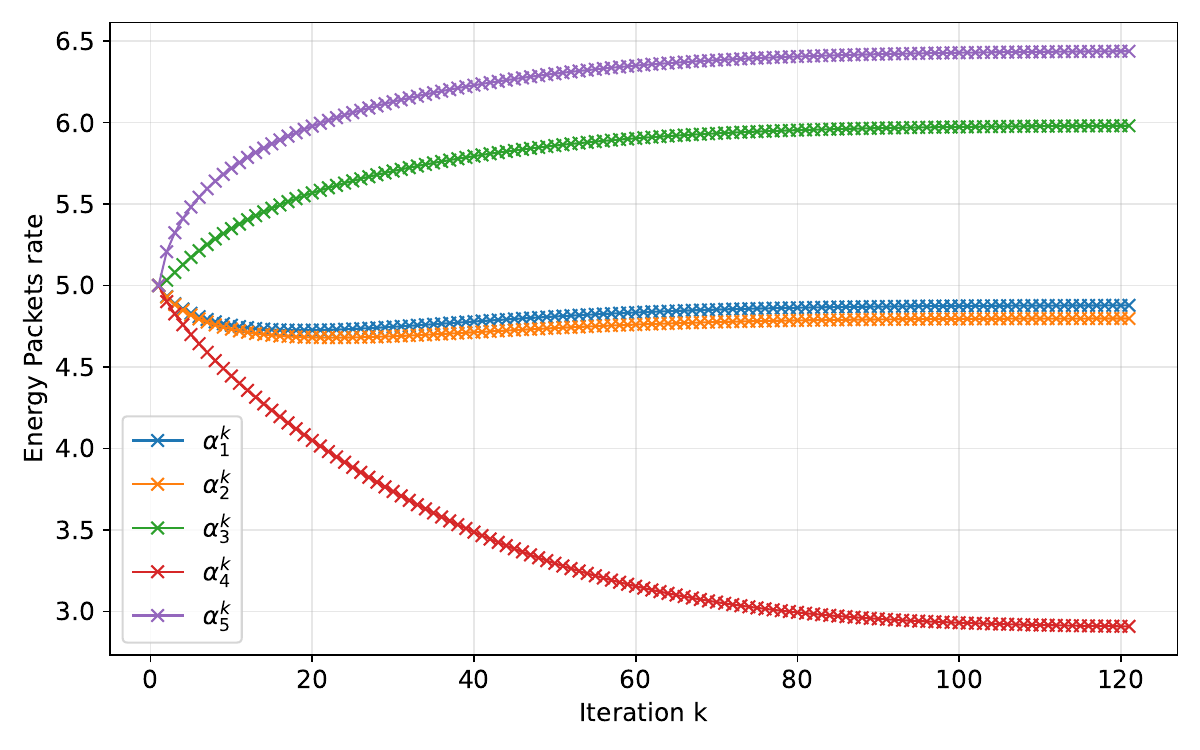}
    \vspace{-0.25cm}
    \caption{Evolution of the energy–allocation parameters 
$\alpha_i^{(k)} = \theta_i^{(k)}$ over the PG--Flow iterations.}
    \label{fig:EPN_alpha}
\end{figure}
\begin{figure}[!ht]
\vspace{-0.5cm}
    \centering
    \includegraphics[width=10cm, height=5.5cm]{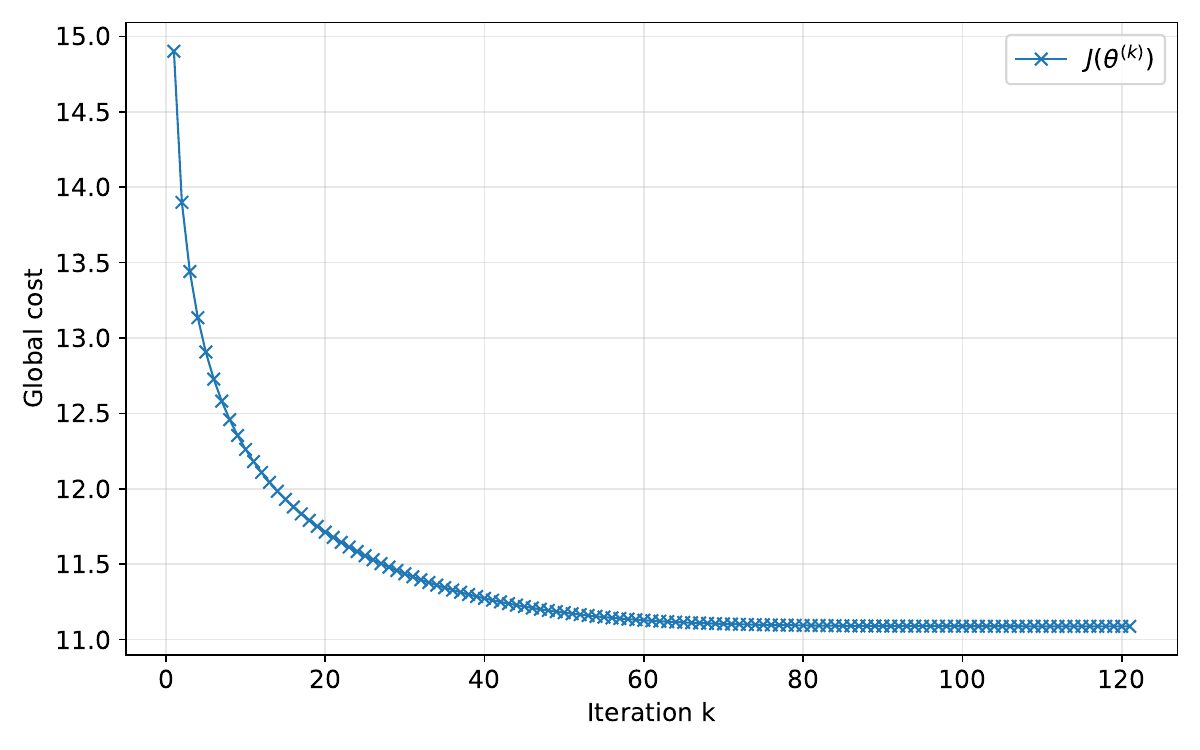}
    \vspace{-0.25cm}
    \caption{Monotone decrease of the objective $J(\theta^{(k)})$ along
PG--Flow iterations.}
    \label{fig:EPN_J}
\end{figure}
\begin{figure}[!ht]
    \vspace{-0.5cm}
    \centering
    \includegraphics[width=11cm, height=5.8cm]{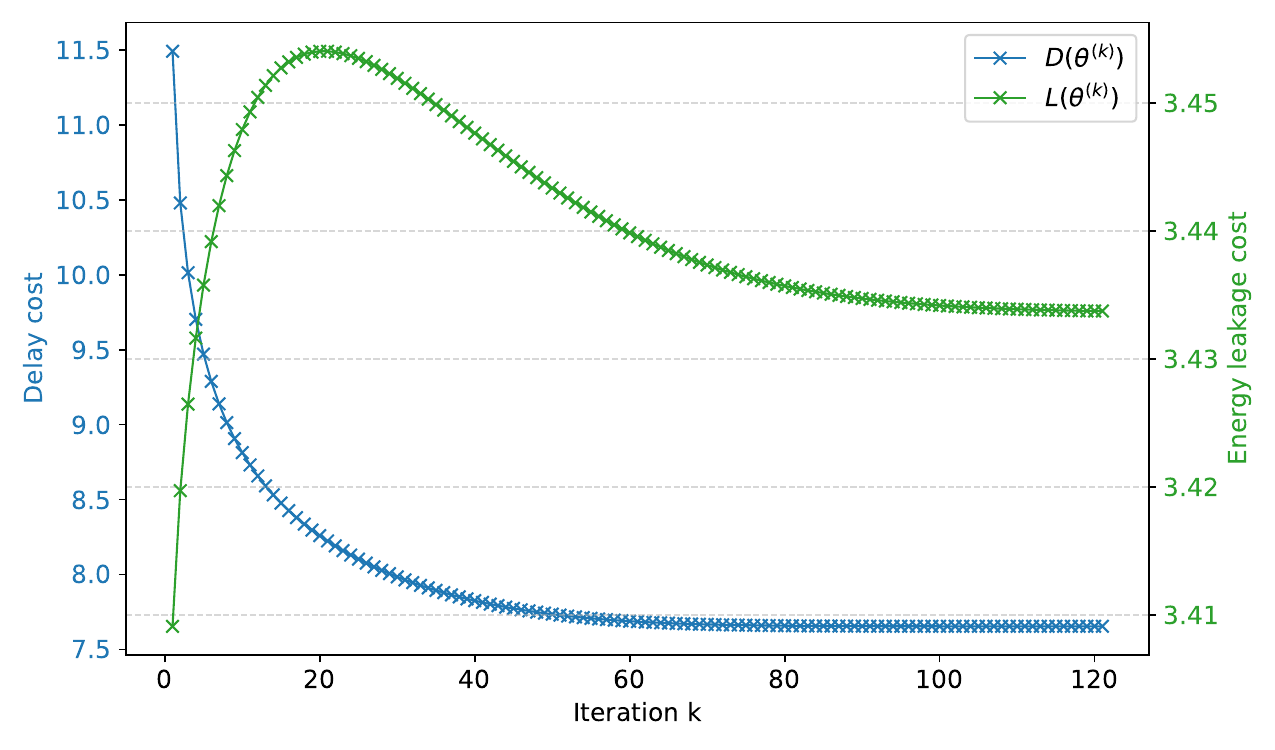}
    \vspace{-0.2cm}
    \caption{Decomposition of the objective into delay 
$D(\theta^{(k)})$ and leakage $L(\theta^{(k)})$ components.}
    \label{fig:EPN_DL}
\end{figure}

Starting from the uniform allocation 
\(
\vect{\theta}^{(0)} = (5,5,5,5,5)
\)
with total budget \(B_{max}= 25\), PG--Flow converges after
about \(120\) iterations to an energy allocation
\[
    \vect{\theta}^\star \approx (4.88,\; 4.80,\; 5.98,\; 2.91,\; 6.44),
    \qquad
    \sum_{i=1}^5 \theta_i^\star = 25.
\]
The corresponding EP loads are
\[
    \vect{\beta}^\star \approx (0.44,\; 0.44,\; 1.00,\; 0.48,\; 1.07),
\]
and the steady-state DP flows remain
\[
    \vect{\phi}^\star \approx (2.64,\; 2.58,\; 3.19,\; 1.27,\; 3.48),
\]
as determined by the Jackson flow equations.
Hence nodes~3 and~5 carry the largest DP loads and are structurally central
in the network (in particular through the cycles involving nodes~1--3--5),
while node~4 carries a comparatively small flow.
It is therefore natural, from the point of view of delay minimization, that
PG--Flow learns to concentrate energy on nodes~3 and~5, to a lesser extent on
nodes~1 and~2, and to allocate very little energy to node~4.  The limiting
allocation \(\vect{\theta}^\star\) can thus be interpreted as a principled
way of ``feeding'' the main bottlenecks of the network under a global budget
constraint.

Figure~\ref{fig:EPN_alpha} illustrates the evolution of the five components
\(\alpha_i^{(k)} = \theta_i^{(k)}\) along the iterations.
Starting from the uniform vector, the parameters quickly separate and then
stabilize, with the clear ordering
\[
    \theta_5^\star \gtrsim \theta_3^\star > \theta_1^\star \approx \theta_2^\star
    \gg \theta_4^\star.
\]
This behavior is fully consistent with the structure of the steady-state DP
flows: the nodes that process the largest amount of traffic (3 and~5) end up
receiving most of the energy, whereas node~4, which sees only a modest flow
and primarily acts as a transient node towards node~5, is assigned a very
small share of the budget.  The trajectories in
Figure~\ref{fig:EPN_alpha} can therefore be read as the progressive
redistribution of energy away from lightly utilized nodes and towards the
main DP bottlenecks. The evolution of the objective \(J(\vect{\theta}^{(k)})\) is reported in
Figure~\ref{fig:EPN_J}.
Initially, we have
\[
    J(\vect{\theta}^{(0)}) 
    = D(\vect{\theta}^{(0)}) + L(\vect{\theta}^{(0)})
    \approx 11.49 + 3.41
    \approx 14.90,
\]
while at the final iterate
\[
    J(\vect{\theta}^\star)
    \approx 11.09
    \quad\text{with}\quad
    D(\vect{\theta}^\star) \approx 7.65,
    \quad
    L(\vect{\theta}^\star) \approx 3.43.
\]
The curve \(k \mapsto J(\vect{\theta}^{(k)})\) decreases monotonically and
shows a typical pattern for gradient-based schemes: a fast decrease during
the first few tens of iterations, followed by a refinement regime in which
the objective approaches its limiting value more slowly.
The small relative improvement observed at the final iterations is
consistent with PG--Flow having reached a point close to a global minimum.

To better understand the underlying trade-off, Figure~\ref{fig:EPN_DL}
displays separately the delay component \(D(\vect{\theta}^{(k)})\) and
the energy leakage \(L(\vect{\theta}^{(k)})\).
At the beginning of the trajectory, \(D\) decreases sharply, while \(L\)
increases moderately: PG--Flow first exploits the available budget more
aggressively, pushing energy towards the heavily loaded nodes in order to
dissipate congestion and reduce queueing delay.
Once the main delay reduction has been achieved, the curve \(D(\vect{\theta}^{(k)})\)
flattens, and subsequent iterations mostly adjust the allocation to re-balance
leakage and delay.
In this later phase, \(L(\vect{\theta}^{(k)})\) slightly decreases again, while
\(D(\vect{\theta}^{(k)})\) remains essentially stable, leading to a final
allocation that realises a sensible compromise between queueing performance
and energy leakage.
Overall, Figures~\ref{fig:EPN_alpha}--\ref{fig:EPN_DL} provide a coherent
behavior: starting from a uniform allocation, PG--Flow reassigns energy to
the critical nodes~3 and~5, yields a substantial reduction of the delay
component, and settles on a stable operating point where the total cost
\(J\) is the closest to its optimal value.

\subsection{Performance evaluation of the algorithms}

In this section, we compare the algorithms introduced in Section~\ref{sec:complexity}. 
(i) \textbf{FD-J} denotes the baseline approach, where $\nabla_{\vect{\theta}} J$ is approximated using finite-difference gradients. 
(ii) \textbf{PG-num} refers to the Policy-Gradient Flow method proposed in this work, but relying exclusively on numerical solvers for the flow fixed point, the adjoint system, and the Jacobian of $\op{G}$. 
(iii) \textbf{PG-analytic} corresponds to the optimized version of PG-Flow specialized to forward (i.e., acyclic) queueing networks, where all derivatives admit closed-form expressions. 
(iv) \textbf{PPO} (Proximal Policy Optimization) represents a simulation-based algorithm relying on Monte Carlo trajectory sampling. PPO  introduced by Schulman et al.~(2017) in \cite{PPO2017}, is a policy–gradient 
method specifically designed to address the stability issues of classical Actor–Critic and 
REINFORCE algorithms. Instead of applying potentially large parameter updates—which often lead to 
high–variance gradients and poor learning stability—PPO constrains successive policies to remain 
close to each other through a clipped surrogate objective. This results in significantly more stable 
and monotonic policy improvement while retaining the simplicity and sample–efficiency of first–order 
methods. We include PPO in our comparison because it is widely regarded as one of the most reliable 
simulation–based policy–gradient algorithms in reinforcement learning: it is substantially more 
robust to hyper-parameters, less sensitive to noise, and typically achieves smoother convergence than 
standard Actor–Critic (A2C/A3C) methods \cite{PPO2017}. For each method, we report, in Table \ref{tab:comp1} and \ref{tab:comp2}, the final global cost $J(\vect{\theta})$, the number of iterations required to reach convergence (or a prescribed precision~$\varepsilon$), and the total execution time (CPU time). 

\textbf{Network generation.}
For the synthetic experiments reported in tables below, we generate random Jackson networks structured as directed
acyclic graphs (DAGs), in the spirit of the construction described in
Subsection~\ref{subsec:jack}. This design enables a fair comparison of all
optimization methods discussed in Section~\ref{sec:complexity}.  

Given the number of queues $d$ and the number of controllable routing links
$p$, we construct a forward DAG in which each node $i$ only routes flow to
nodes $j>i$. External arrivals occur only at node $0$, with rate
$\lambda^{\mathrm{ext}}_0 = 4$. Service rates alternate between ``slow'' and
``fast'' servers, i.e., $\mu_i \in \{8,12\}$ according to the parity of $i$. Among the nodes $\{1,\dots,d-2\}$, a subset of size $p$ is randomly selected as controllable branching nodes. At each such node $i$, we create two
outgoing arcs: a ``slow'' arc to $i+1$ and a ``fast'' arc to a randomly chosen
successor $j>i+1$. Their routing probabilities are parametrized by a control
variable $\theta_k$, namely $P(i\to i+1)=1-\theta_k$ and $P(i\to j)=\theta_k$.
All remaining nodes are non-controllable and route flow according to fixed
probabilities: a random value $x\in(0.2,0.8)$ determines $P(i\to i+1)=x$ and
$P(i\to j)=1-x$ for a randomly selected $j>i+1$. Node $d-1$ has a single fixed
successor $d$, and node $d$ has no outgoing arcs, meaning that all departures
from $d$ immediately leave the network. 

\textbf{Openess.}
This construction guarantees that the network is open for all
$\vect{\theta}\in\op{U}$ (with $\op{U} = [0,1]^p$). That is, for every $i<d$, every queue
admits a directed path $i \to\ i+1 \to i+2 \to \dots \to\ d$
with strictly positive probability. Thus every job exits the system in a finite number of steps. Finally, the presence of fixed (non-controllable) routing links plays an important role: they act as \emph{structural constraints} that limit the action space and make the resulting optimization problem more realistic and challenging.

\textbf{Stability.}
In all generated instances, stability is guaranteed for any 
$\vect{\theta}\in\op{U}$ . Since external arrivals occur only at node~1 with 
rate $4$ and the routing graph is a forward DAG. Hence the effective 
arrival rate satisfies $\phi_i^\star(\vect{\theta}) \le 4$ for all $i$. With 
service rates $\mu_i\in\{8,12\}$, we have $\rho_i(\vect{\theta}) \;=\; \frac{\phi_i^\star(\vect{\theta})}{\mu_i}
\;\le\; \frac{4}{8}  \;<\; 1$. 
Hence the network is stable for all admissible 
parameters. 

Other hypothesis verifications (H3) and (H4) follow the same arguments as in Appendix~\ref{annexe:A}.

\textbf{Objective.}
We fixed the optimized quantity as the total mean number of tasks in 
the system, i.e., the standard Jackson-network cost
$J(\vect{\theta})$ defined in~\eqref{eq:jackson-2param-cost}.

\textbf{Hyper-parameter setting.} 
We fixed the same stopping condition for all optimization methods (FD--J, PG--num and PG--analytic). At iteration $k$, the algorithm stops whenever
\[
\frac{|J(\vect{\theta}_k)-J(\vect{\theta}_{k-1})|}{\max(1,|J(\vect{\theta}_{k-1})|)} 
\leq \epsilon_J
\qquad\text{or}\qquad
\|\nabla J(\vect{\theta}_k)\|_2 \leq \epsilon_{\nabla},
\]
where the first test monitors the \emph{relative} improvement of the objective, 
and the second corresponds to a first-order optimality condition.  
The tolerances were chosen so as to ensure reliable convergence without 
performing unnecessary iterations: we use 
$\epsilon_J = 10^{-6}$, $\epsilon_{\nabla} = 10^{-4}$, 
and $\texttt{max\_Iter} = 500$. Additional internal tolerances are required for implicit-differentiation 
methods. In PG--num, the steady-state flows 
$\vect{\phi}^\star(\vect{\theta})$ are computed by Anderson acceleration with a 
fixed-point tolerance $\epsilon_{\vect{\phi}} = 10^{-10}$, and all derivatives 
$\partial_{\vect{\phi}}$ and $\partial_{\vect{\theta}}$ are approximated by finite differences using a step 
$\epsilon_{\partial} = 10^{-8}$.  
For FD--J, the finite-difference approximation is applied directly to the full
objective $J$, using the same internal finite-difference step $\epsilon_{\partial}$. For all methods, the parameters are updated through a projected gradient step
with a constant step-size $\eta = 0.05$. 

For the Monte-Carlo simulation based baseline (PPO), the hyper-parameters must 
be chosen with care to obtain stable gradients and a reasonable variance–bias 
trade-off. In all experiments, we fix the number of policy updates to 
$n_{\text{iter}}=50$, and each update is based on a batch of 
$N_{\text{traj}}=100$ simulated trajectories. Each trajectory is generated over a 
finite time horizon $H \in\{50, 100\}$ (refer to Table \ref{tab:comp1}). For each 
batch of trajectories, the PPO objective is optimized for $n_{\text{epochs}} = 4$ 
epochs using the Adam optimizer with learning rate $10^{-2}$ and 
clipping parameter $\varepsilon_{\text{clip}} = 0.2$. These values provided a 
good balance between computational cost and stability, while avoiding the high 
variance issues typically observed when using longer horizons or smaller batches. 

\begin{table}[h]
  \centering
  \caption{Performance comparison for small and medium-scale Jackson DAGs. Where $d$ the number of queues and $p$ the number of control parameters. \label{tab:comp1}}
  \begin{tabular}{cc l || r r r}
    $d$ & $p$ & Method 
        & $J(\vect{\theta^*})$ 
        & Iterations 
        & Exec. time (s) \\
    \midrule
    \multirow{4}{*}{10} & \multirow{4}{*}{3} 
    & PG-analytic    &  $2.4737$ & $106$ & $0.006$ \\
    & & PG-num       &  $2.4737$ & $106$ & $0.396$ \\
    & & FD-J         &  $2.4737$ & $106$ & $0.263$ \\
    & & PPO ($H$=50)   &  $2.5150$ & $50$ & $42.24$ \\
    \midrule
    \multirow{4}{*}{50} & \multirow{4}{*}{20}
    &  PG-analytic     &  $5.3015$ & $346$ & $0.188$ \\
    & & PG-num         &  $5.5722$ & $500$ & $15.66$ \\
    & & FD-J           &  $5.5811$ & $500$ & $33.49$ \\
    & & PPO ($H=$50)     &  $5.4766$ & $50$ & $117.19$ \\
    \midrule
    \multirow{4}{*}{100} & \multirow{4}{*}{40}
    & PG-analytic     &  $4.9482$ & $263$ & $0.237$ \\
    & & PG-num        &  $5.4432$ & $500$ & $160.50$ \\
    & & FD-J          &  $5.0723$ & $500$ & $1399.61$ \\
    & & PPO ($H$=100)   &  $5.2886$ & $50$ & $174.85$ \\
    \bottomrule
  \end{tabular}
\end{table}

\begin{table}[h]
\centering
\caption{
Scalability of the PG-analytic method on large Jackson DAGs. Where $d$ the number of queues and $p$ the number of control parameters. \label{tab:comp2}}
\label{tab:scaling_analytic}
\begin{tabular}{c|c||c c c}
 $d$ & $p$ & $J(\vect{\theta^*})$  & Iterations & Exec. time (s) \\
\hline
\multirow{2}{*}{$5\,000$ } 
& $2000$   & $5.64$ & $707$ & $29.10$ \\
& $4000$   & $4.69$ & $202$ & $9.28$ \\ \hline
\multirow{2}{*}{$10\,000$ } 
& $3000$ & $6.96$ & $1388$ & $172.82$ \\
& $9000$ & $4.73$ & $724$ & $95.55$ \\ \hline

\multirow{2}{*}{$50\,000$ } 
& $20\,000$ & $7.07$ & $1318$ & $1252.43$ \\ 
& $40\,000$ & $5.32$ & $130$ & $262.05$ \\ \hline

\multirow{2}{*}{$100\,000$ } 
 & $50\,000$ & $7.32$ & $1330$ & $1022.68$ \\
 & $80\,000$ & $4.98$ & $111$ & $429.79$ \\
\hline

$200\,000$ & $150\,000$ & $5.91$ & $355$ & $1739.83$ \\ \hline
\end{tabular}
\end{table}

\textbf{Analysis of results.} 
Table~\ref{tab:comp1} compares the four optimization approaches on small and medium-sized Jackson DAGs. The PG-analytic method provides a useful reference point, as it computes exact gradients with a linear-time complexity in acyclic networks (Equation~\eqref{eq:O-analytic}). Its execution time remains very low across all problem sizes, and it consistently satisfies the stopping criterion. The PG-num and FD--J methods also produce accurate solutions, especially for small networks, but their reliance on finite-difference perturbations introduces numerical sensitivity as $d$ grows. This effect is visible in the objective values for $d=50$ and $d=100$, where both methods reach the maximum iteration (500 iterations) without meeting the targeted precision $\epsilon_J = 10^{-5}$. Their execution times align with the theoretical complexities discussed in Section~\ref{sec:complexity}, which increases computational cost for larger dimensions. The PPO method produces approximate solutions across all settings and returns objective values of the same general order as the finite-difference approaches. Its performance depends on several hyper-parameters, such as the horizon length $H$ and the number of sampled trajectories $N_{\text{traj}}$, and each iteration requires evaluating approximately $N_{\text{traj}} \times H$ simulated transitions. For the small and medium networks considered here, the resulting execution time remains reasonable, but PPO becomes substantially slower in larger systems, as longer horizons are typically required for trajectories to reach the steady-state behavior of the queues. For instance, a simulation of $d=500$ and $p=300$ yields $4505.8$ seconds (with $H=200$ and $N_{traj}=100$) with $J(\theta)=4.56$, while PG-analytic converges in $0.45$ seconds achieving $J(\theta)=4.52$. Overall, PPO provides meaningful approximate solutions, while PG-analytic achieves the most favorable balance between accuracy, convergence, and computational efficiency.

Given the promising results of PG-analytic in Table~\ref{tab:comp1}, we further investigate its scalability in Table~\ref{tab:comp2}. The execution times grow approximately linearly with~$d$, in accordance with the per-iteration complexity of Equation~\eqref{eq:O-analytic} in the acyclic setting. This behavior confirms that PG-analytic remains suitable for problems involving tens of thousands of queues and beyond.

A second observation concerns the influence of the number of control parameters~$p$. When $p$ is relatively small (compared to $d$), such as $d = 100{,}000$ and $p = 50{,}000$ (or even lower), each iteration is computationally inexpensive, but the optimization typically requires more iterations to reach the stopping criterion. This behavior is expected: fixing a large subset of non-controlled routing parameters reduces the dimensionality of the control space and therefore limits the flexibility of the optimization process. With fewer degrees of possibilities, additional iterations are needed to compensate for these structural constraints.

Conversely, when $p$ is larger—for example $p = 80{,}000$ for $d = 100{,}000$—the optimization benefits from a richer parametrization and converges in substantially fewer iterations, despite a slightly higher per-iteration computation time. A related observation is that the objective value $J(\vect{\theta}^*)$ decreases as $p$ increases. This is consistent with the fact that enlarging the controllable parameter space provides more opportunities to reduce the total cost, whereas fixed (non-controllable) parameters effectively act as constraints. Overall, these results confirm both the scalability and the practical robustness of the PG-analytic approach, while highlighting the important interplay between network size and policy parametrization in large-scale optimization. 

Experiments were conducted on a laptop equipped with 10 CPU cores (8 cores at 3.2~GHz peak frequency and 2 cores at 2.0~GHz), and 16~GB of RAM.  

\section{Conclusion}
\label{sec:conclusion}

We introduced PG--Flow, a deterministic policy-gradient framework for steady-state
optimization in geometric product-form queueing networks.  
By expressing classical performance metrics in terms of local steady-state flows,
the method leverages the fixed-point structure of product-form models and computes
exact policy gradients via implicit differentiation and a local adjoint system.
We established global convergence under affine flow operators and natural convexity
assumptions, and identified acyclic network structures that enable linear-time
implementations.  
Numerical experiments on routing control in Jackson networks and on energy-arrival
control in Energy Packet Networks demonstrated the effectiveness and scalability of
the approach.

Several directions for future work arise naturally.  
One avenue is to extend PG--Flow beyond geometric product-form networks, either by
considering richer steady-state parametrizations or by incorporating approximate
flow representations.  
Another promising direction is the integration of PG--Flow within adaptive or
online control architectures, combining deterministic gradients with real-time
measurements.  
Finally, exploiting sparsity, decomposability, or hierarchical structures may allow further reductions in computational cost.  
We hope that PG--Flow contributes to strengthening the interface between analytical
queueing theory and modern gradient-based methods for steady-state control.


\bibliography{PG-Flow}

\appendix
\section{Appendix: Jackson Network \label{annexe:A}}
\subsubsection*{Verification of the PG--Flow assumptions}

We verify that this network satisfies the structural conditions of
Section~\ref{sec:global-convergence}.

\paragraph{(H1) Convex and compact domain.}
We consider the rectangular parameter domain
\[
\mathcal U
= [\theta_{1,\min},\theta_{1,\max}] \times
  [\theta_{2,\min},\theta_{2,\max}]
= [0,1]\times[0,1],
\]
which is nonempty, convex, and compact.
The routing probabilities are
$p(\theta_1)=\theta_1$, $q(\theta_2)=\theta_2,$
and therefore lie in $[0,1]$ for all $\vect{\theta}\in \op{U}$.
The steady-state flows are
\begin{equation}
\label{eq:num:closed}
\Lambda_1 = \lambda_1^{\mathrm{ext}},\qquad
\Lambda_2 = \lambda_1^{\mathrm{ext}}\,p(\theta_1),\qquad
\Lambda_3 = \lambda_1^{\mathrm{ext}}\bigl(1-p(\theta_1) + q(\theta_2)p(\theta_1)\bigr),
\end{equation}

which satisfy
\[
0\le \Lambda_2\le \lambda_1^{\mathrm{ext}},\qquad
0\le \Lambda_3\le \lambda_1^{\mathrm{ext}}.
\]
Since $\lambda_1^{\mathrm{ext}} = 4$ then  $\Lambda_i<\mu_i=(6,5,7)$ for all $\vect{\theta}\in \op{U}$, the network is stable throughout $\op{U}$.
Hence (H1) holds.

\paragraph{(H2) Openness.}
By definition, the dependency network $A(\vect{\theta})$  is feed-forward, hence for every queue there is a path to a departure queue. This structure holds for every $\vect{\theta}\in \op{U}$. 

\paragraph{(H3) Monotonicity of the steady-state flows.}
In this acyclic Jackson network, the steady-state flows admit the closed-form expressions \eqref{eq:num:closed} on $\op{U}$, the routing probabilities
$p(\theta_1)=\theta_1$ $q(\theta_2)=\theta_2$
are coordinate-wise non-decreasing.
Differentiating the steady-state flows yields
\[
\frac{\partial \Lambda_2}{\partial p} = 4 \ge 0,\qquad
\frac{\partial \Lambda_3}{\partial p} = 4\,q(\theta_2)\ge 0,\qquad
\frac{\partial \Lambda_3}{\partial q} = 4\,p(\theta_1)\ge 0.
\]
Hence, for any $\vect{\theta}\le\vect{\theta'}$ (component-wise), $\Lambda_i(\vect{\theta})\le \Lambda_i(\vect{\theta}')$ for $i\in\{1,2,3\}$.
Thus the steady-state flow vector $\vect{\phi}^\star(\vect{\theta})$ is coordinate-wise non-decreasing in the
control parameters, and (H3) holds.

\paragraph{(H4) Convex and monotone local rewards.}
Each queue $i\in\{1,2,3\}$ has local reward
$r_i(x)=\frac{x}{\mu_i-x},\qquad x\in[0,\mu_i),$
which is convex and strictly increasing.
The steady-state objective can be written as
$J(\vect{\theta})=\sum_{i=1}^3 r_i(\Lambda_i^\star(\vect{\theta}))$.
Since each $r_i$ is convex and non-decreasing and each $\Lambda_i^\star(\vect{\theta})$ is
coordinate-wise non-decreasing by (H3), the composition $r_i(\Lambda_i^\star(\vect{\theta}))$ is
convex on $\op{U}$.
Therefore (H4) holds.

\medskip
In conclusion, all assumptions (H1)--(H4) of Theorem~\ref{th:global-PGFlow} are satisfied
for this two-parameter Jackson network. In the next subsection, we provide details of the first iteration of PG-flow.

\subsubsection*{First PG-Flow iteration}

We initiate PG--flow control parameters with
\(\vect{\theta}_0 = (0.8,0.8)\).
hence
\[
p_0 = p(\theta_{1,0}) = 0.8,\qquad
q_0 = q(\theta_{2,0}) = 0.8.
\]

\paragraph{1. Flow equations.}
At \(\theta_0\), the steady-state flows are
\[
\Lambda_1^{(0)} = 4,\qquad
\Lambda_2^{(0)} = 4p_0 = 3.2,\qquad
\Lambda_3^{(0)} = 4\bigl(1-p_0+q_0p_0\bigr)
= 3.36.
\]
Thus
\[
\boxed{
\vect{\phi}^{(0)}
=
\begin{bmatrix}
4\\[3pt]3.2\\[3pt]3.36
\end{bmatrix}}.
\]

\paragraph{2. Objective value.}
The initial cost is
\[
J(\vect{\theta}_0)
=
\frac{4}{6-4}
+\frac{3.2}{5-3.2}
+\frac{3.36}{7-3.36}
\approx \boxed{4.70}.
\]

\paragraph{3. Local Jacobians.}
The flow operator is described in \eqref{eq:num:fixed}, hence the Jacobian with respect to $\vect{\phi}$ is
\[
G_{\vect{\phi}}(\vect{\theta})
=
\frac{\partial \op{G}}{\partial \vect{\phi}}
=
\begin{bmatrix}
0 & 0 & 0\\[3pt]
p(\theta_1) & 0 & 0\\[3pt]
1-p(\theta_1) & q(\theta_2) & 0
\end{bmatrix}.
\]
At \(\vect{\theta}_0\) (i.e., \(p_0=q_0=0.8\)):
\[
\boxed{
G_{\vect{\phi}}({\vect{\theta}}_0)
=
\begin{bmatrix}
0 & 0 & 0\\[3pt]
0.8 & 0 & 0\\[3pt]
0.2 & 0.8 & 0
\end{bmatrix}},
\qquad
I - G_{\vect{\phi}}({\vect{\theta}}_0)
=
\begin{bmatrix}
1 & 0 & 0\\[3pt]
-0.8 & 1 & 0\\[3pt]
-0.2 & -0.8 & 1
\end{bmatrix}.
\]

The Jacobian w.r.t.~\(\vect{\theta}\) is obtained from
\[
\frac{\partial \op{G}}{\partial \theta_1}
=
\begin{bmatrix}
0\\[3pt]
(\partial_{\theta_1}p)\,\phi_1\\[3pt]
-(\partial_{\theta_1}p)\,\phi_1
\end{bmatrix},
\qquad
\frac{\partial \op{G}}{\partial \theta_2}
=
\begin{bmatrix}
0\\[3pt]
0\\[3pt]
(\partial_{\theta_2}q)\,\phi_2
\end{bmatrix}.
\]
Since \(p(\theta_1)=\theta_1\) and \(q(\theta_2)=\theta_2\),
we have \(\partial_{\theta_1}p = 1\), \(\partial_{\theta_2}q = 1\).
Using \(\vect{\phi}^{(0)}=(4,3.2,3.36)\), we obtain
\[
\boxed{
G_{\vect{\theta}}(\vect{\theta}_0)
=
\begin{bmatrix}
0 & 0\\[2pt]
4 & 0\\[2pt]
-4 & 3.2
\end{bmatrix}}
\]

\paragraph{4. Objective derivatives.}
The objective \(\op{F}(\vect{\phi})=\sum_i \tfrac{\Lambda_i}{\mu_i-\Lambda_i}\) yields
\[
\frac{\partial \op{F}}{\partial \Lambda_i}
=
\frac{\mu_i}{(\mu_i-\Lambda_i)^2}.
\]
Evaluated at \(\vect{\phi}^{(0)}=(4,3.2,3.36)\):
\[
\frac{\partial \op{F}}{\partial \Lambda_1}
= \frac{6}{(6-4)^2} = 1.5, \ \ \ \
\frac{\partial \op{F}}{\partial \Lambda_2}
= \frac{5}{(5-3.2)^2} \approx 1.5432, \ \ \ \
\frac{\partial \op{F}}{\partial \Lambda_3}
= \frac{7}{(7-3.36)^2} \approx 0.5283.
\]
Thus
\[
\boxed{
\partial_{\vect{\phi}} \op{F}(\vect{\phi}^{(0)})
\approx
\begin{bmatrix}
1.50 & 1.54 & 0.53
\end{bmatrix}},
\qquad
\partial_{\vect{\theta}} \op{F}(\vect{\phi}^{(0)},\vect{\theta}_0)
=
\begin{bmatrix}
0 & 0
\end{bmatrix},
\]
since \(\op{F}\) does not depend explicitly on~\(\vect{\theta}\).

\paragraph{5. Adjoint system.}
The adjoint variable \(y\in\mathbb R^3\) solves $(I - G_{\vect{\phi}}(\vect{\theta}_0))^\top y = (\partial_{\vect{\phi}} \op{F})^\top.$

We have
\[
(I - G_{\vect{\phi}}(\vect{\theta}_0))^\top =
\begin{bmatrix}
1 & -0.8 & -0.2\\
0 & 1 & -0.8\\
0 & 0 & 1
\end{bmatrix},
\qquad
(\partial_{\vect{\phi}} \op{F})^\top \approx
\begin{bmatrix}
1.50\\[2pt]1.54\\[2pt]0.53
\end{bmatrix}.
\]
Solving this upper-triangular system (back-substitution) yields
\[
y_3 \approx 0.53,\qquad
y_2 \approx 1.54 + 0.8\,y_3 \approx 1.97, \qquad
y_1 \approx 1.50 + 0.8\,y_2 + 0.2\,y_3
\approx 3.18.
\]
Hence
\[
\boxed{
y \approx
\begin{bmatrix}
3.18\\[2pt]
1.97\\[2pt]
0.53
\end{bmatrix}}.
\]

\paragraph{6. Global gradient.}
The implicit gradient is
\[
\nabla_{\vect{\theta}} J(\vect{\theta}_0)
= \partial_{\vect{\theta}} \op{F} + G_{\vect{\theta}}(\vect{\theta}_0)^\top y
= G_{\vect{\theta}}(\vect{\theta}_0)^\top y.
\]
Using
\[
G_{\vect{\theta}}(\vect{\theta}_0)^\top y
=
\begin{bmatrix}
0 & 4 & -4\\[2pt]
0 & 0 & 3.2
\end{bmatrix},
\]
we obtain
\[
\nabla_{\theta_1} J(\vect{\theta}_0)
= 4y_2 -4y_3
\approx 5.75,
\]
\[
\nabla_{\theta_2} J(\vect{\theta}_0)
= 3.2\,y_3
\approx 1.69.
\]
Thus
\[
\boxed{
\nabla_{\vect{\theta}} J(\vect{\theta}_0)
\approx
\begin{bmatrix}
5.75 & 1.69
\end{bmatrix}}.
\]

\paragraph{7. Update and new cost.}
With a step-size \(\eta = 0.01\), the PG-Flow update is
\[
\vect{\theta_1}
= \Pi_{\mathcal{U}}\Big( \vect{\theta_0} - \eta\,\nabla_{\vect{\theta}} J(\vect{\theta}_0) \Big)
= \Pi_{\mathcal{U}}\Big(
\begin{bmatrix}
0.8\\[2pt]0.8
\end{bmatrix}
- 0.01
\begin{bmatrix}
5.7502\\[2pt]1.6906
\end{bmatrix}
\Big)
\approx
\begin{bmatrix}
0.7425\\[2pt]0.7831
\end{bmatrix}.
\]
The corresponding routing probabilities are
\[
p_1 = p(\vect{\theta}_{1,1}) \approx 0.7425,\qquad
q_1 = q(\vect{\theta}_{2,1}) \approx 0.7831.
\]
The new steady-state flows are
\[
\Lambda_1^{(1)} = 4,
\qquad
\Lambda_2^{(1)} = 4p_1 \approx 2.97,
\qquad
\Lambda_3^{(1)} = 4\bigl(1-p_1 + q_1 p_1\bigr)
\approx 3.36.
\]
The updated cost is
\[
J(\vect{\theta}_1)
=
\frac{4}{6-4}
+\frac{\Lambda_2^{(1)}}{5-\Lambda_2^{(1)}}
+\frac{\Lambda_3^{(1)}}{7-\Lambda_3^{(1)}}
\approx 4.38.
\]
Thus
\[
J(\vect{\theta}_1) \approx 4.38 \;<\; J(\vect{\theta}_0)\approx 4.70,
\]
which shows a strict decrease of the objective after a single PG-Flow
iteration, in line with the monotone descent guarantee of
Theorem~\ref{th:global-PGFlow}.

\section{Appendix: EPN Network}
\label{annexe:B}

\subsubsection*{Verification of the PG--Flow assumptions.}

We verify assumptions (H1)--(H4) from Section~\ref{sec:global-convergence}.

\begin{itemize}
        \item[(H1)] 
    The feasible control domain is a nonempty compact convex set.
    In the EPN model, the control vector 
    $\alpha(\vect{\theta})=(\alpha_1(\vect{\theta}),\ldots,\alpha_N(\vect{\theta}))$
    represents energy allocation intensities and $B_{\max}$ is the maximum energy budget to be allocated across the network.
    Hence 
    \[
        \alpha_i(\vect{\theta})\;\ge 0,
        \qquad 
        \sum_{i=1}^N \alpha_i(\vect{\theta})\;\le B_{\max},
    \]
defines a closed, bounded, and convex simplex under a budget constraint. 
In particular, the Euclidean projection onto this feasible set is well defined.
    \item[(H2)] 
    Since routing $P$ is fixed and forms an open Jackson network, the DP flow system
    admits a unique solution.

    \item[(H3)] 
In the EPN model, the steady-state DP flows $\phi_i^\star$ satisfy the
linear Jackson system $\phi^\star = \lambda^{\mathrm{ext}} + P^\top \phi^\star,$ whose solution does not depend on $\vect{\theta}$, since routing is fixed.
Hence $\phi^\star(\vect{\theta})$ is constant on $\op{U}$ and therefore
coordinate-wise non-decreasing.  
Moreover, the EP flows are given by
$\beta_i(\vect{\theta})=\frac{\alpha_i(\vect{\theta})}{\gamma_i+\mu_i},$
and since $\alpha_i(\vect{\theta})$ is affine and non-decreasing in each
component of $\vect{\theta}$, so is $\beta_i(\vect{\theta})$.  
Thus the steady-state flow vector
\[
\vect{\phi}^\star(\vect{\theta}) = (\phi_1^\star,\dots,\phi_N^\star,\;
\beta_1(\vect{\theta}),\dots,\beta_N(\vect{\theta}))
\]
is coordinate-wise non-decreasing in $\vect{\theta}$, and (H3) holds.

\item[(H4)] 
For each node $i$, the local cost
$r_i(\phi_i,\beta_i)=\frac{\phi_i}{\mu_i\beta_i-\phi_i}$
is convex and non-increasing in $\beta_i$ for fixed $\phi_i$, while the
leakage term $\gamma_i\beta_i(\vect{\theta})$ is linear.  
In the present EPN model, the DP flows $\phi_i^\star$ are independent of 
$\vect{\theta}$ and the EP flows $\beta_i(\vect{\theta})$ are affine in
$\vect{\theta}$.  
Hence the composed map 
\[
\vect{\theta} \mapsto r_i(\phi_i^\star,\beta_i(\vect{\theta}))
\]
is convex in $\vect{\theta}$, 
\emph{even though $r_i(\phi_i,\beta_i)$ is not jointly convex in 
$(\phi_i,\beta_i)$ as stipulated in (H4)}.  
Summing over all nodes yields a convex objective $J(\vect{\theta})$,
so the conclusion of (H4) holds for this EPN example.

\end{itemize}

Since (H1)--(H4) are satisfied, PG--Flow applies directly and produces globally convergent
updates of the energy input parameters $\alpha_i(\vect{\theta})$ in this EPN.

\subsubsection*{Numerical Example: 5-node Energy Packet Network}

We now consider a concrete EPN composed of $N=5$ nodes.
Data packets (DPs) circulate according to a fixed routing matrix $P$, while
Energy packets arrivals are controlled through $\alpha_i(\vect{\theta})=\theta_i$.
The global constraint $\sum_{i=1}^N \alpha_i(\vect{\theta})\le B_{\max}$ 
recalls that a finite energy budget must be allocated across the network.
This example illustrates PG--Flow on a topology containing multiple cycles.

\smallskip
\textbf{\textit{DP dynamics.}}
External arrivals occur at rates
\[
    \lambda^{\mathrm{ext}} = (2.0,\; 1.0,\; 0.5,\; 0.5,\; 1.0).
\]
Each DP served from node $i$ is sent to node $j$ with probability $P_{i,j}$, or leaves the network.
The routing probabilities are:
\begin{itemize}
    \item Node 1: $P_{1,2}=0.6$, exit with prob.\ $0.4$;
    \item Node 2: $P_{2,3}=0.5$, $P_{2,4}=0.3$, exit with prob.\ $0.2$;
    \item Node 3: $P_{3,1}=0.2$, $P_{3,5}=0.5$, exit with prob.\ $0.3$;
    \item Node 4: $P_{4,5}=0.7$, exit with prob.\ $0.3$;
    \item Node 5: $P_{5,3}=0.4$, exit with prob.\ $0.6$.
\end{itemize}
This corresponds to the routing matrix:
\[
P = 
\begin{pmatrix}
0 & 0.6 & 0   & 0   & 0   \\
0 & 0   & 0.5 & 0.3 & 0   \\
0.2 & 0 & 0   & 0   & 0.5 \\
0 & 0   & 0   & 0   & 0.7 \\
0 & 0   & 0.4 & 0   & 0
\end{pmatrix}.
\]

\smallskip
\textbf{\textit{EP dynamics}}
Maximum energy budget is fixed to $25$ EPs. We take:
\[
\mu = (10,10,5,5,5),
\qquad
\gamma_i = 1 \;\; \forall i.
\]
Thus
\[
    \beta_i(\vect{\theta}) = \frac{\theta_i}{\gamma_i+\mu_i}
    = \frac{\theta_i}{\mu_i+1}.
\]

We initialise at
\[
    \vect{\theta}^{(0)}=(5,5,5,5,5),
\qquad
    \vect{\beta}^{(0)} \approx (0.4545,\;0.4545,\;0.8333,\;0.8333,\;0.8333).
\]
Hence the effective DP service rate is
\[
    \mu_i\beta_i^{(0)}
    \approx
    (4.545,\;4.545,\;4.167,\;4.167,\;4.167).
\]
In the following we give essential steps of the first iteration of PG-Flow. 

\paragraph{1. DP flow equations}

The DP flows satisfy the Jackson system
$(I - P^\top)\vect{\phi} = \lambda^{\mathrm{ext}}.$
Contrary to an acyclic Jackson network, where the DP flows can be computed in linear algorithmic complexity, the presence of cycles in the routing matrix prevents such a reduction.  
Therefore, the steady-state flows are obtained with an iterative fixed-point method,
whose per-iteration cost is $O(d^2)$ (and whose overall complexity depends on the
number of iterations required). One may also use a direct dense linear solver
(e.g., Gaussian elimination) with a cubic cost $O(d^3)$.
\[
\vect{\phi}^\star \approx (2.64,\; 2.58,\; 3.19,\; 1.27,\; 3.48).
\]
And the stability holds since $\mu_i\beta_i^{(0)} > \phi_i^\star$ for all $i$.

\paragraph{2. Performance at the initial point}

The mean DP queue lengths are
\[
D_i(\vect{\theta})
    =
    \frac{\phi_i^\star}{\mu_i\beta_i(\theta)-\phi_i^\star},
\]
giving
\[
D(\vect{\theta}^{(0)})
    = \sum_i D_i(\theta^{(0)})
    \approx 11.49.
\]

The EP leakage cost is
\[
L(\vect{\theta}^{(0)})
=
\sum_i \gamma_i\beta_i^{(0)}
\approx 3.41.
\]

With $w_1=w_2=1$,
\[
J(\vect{\theta}^{(0)})
= D(\vect{\theta}^{(0)}) + L(\vect{\theta}^{(0)})
\approx 14.90.
\]

\paragraph{3. Gradient and adjoint}

We compute
\[
\frac{\partial \op{F}}{\partial \phi_i}
=
\frac{\mu_i\beta_i}{(\mu_i\beta_i-\phi_i)^2},
\qquad
\frac{\partial \op{F}}{\partial \beta_i}
=
-\frac{\mu_i\phi_i}{(\mu_i\beta_i-\phi_i)^2}+\gamma_i.
\]

At $(\vect{\phi}^\star,\beta^{(0)})$:
\[
\frac{\partial \op{F}}{\partial \vect{\phi}}
\approx
(1.25,\;1.18,\;4.37,\;0.50,\;8.84),
\]
\[
\frac{\partial \op{F}}{\partial \vect{\beta}}
\approx
(-6.27,\;-5.68,\;-15.72,\;0.24,\;-35.90).
\]

The PG--Flow adjoint solves
\[
(I-P) y^{\phi}
=
\left(\frac{\partial F}{\partial \phi}\right)^\top.
\]
This gives
\[
y^{\phi}
\approx
(7.68,\; 10.72,\; 12.90,\; 10.30,\; 13.99),
\qquad
y^{\beta} 
=
\left(\frac{\partial F}{\partial \beta}\right)^\top.
\]

The gradient is
\[
\frac{\partial J}{\partial \theta_i}
=
\frac{y^\beta_i}{\mu_i+1},
\]
yielding
\[
\nabla_{\vect{\theta}} J(\vect{\theta}^{(0)})
\approx
(-0.57,\; -0.52,\; -2.62,\; 0.04,\; -5.98).
\]

\paragraph{4. Update and new cost}

With step-size $\eta=0.05$:
\[
\vect{\theta}^{(1)}
=
\Pi_{\op{U}}\Big( \vect{\theta}^{(0)} - \eta\,\nabla J(\vect{\theta}^{(0)})\Big)
\approx
(4.931,\;4.928,\;5.032,\;4.900,\;5.206).
\]

Re-evaluating the cost:
\[
D(\vect{\theta}^{(1)}) \approx 10.47,
\qquad
L(\vect{\theta}^{(1)}) \approx 3.42,
\]
\[
J(\vect{\theta}^{(1)}) \approx 13.89
< 
J(\vect{\theta}^{(0)}) \approx 14.90.
\]

\end{document}